\documentclass[reprint,onecolumn]{article}
\usepackage[latin9]{inputenc}   
\usepackage{mathtools}
\usepackage{amstext}
\usepackage{amsthm}
\usepackage{graphicx}
\usepackage{dcolumn}
\usepackage{amsfonts} 
\usepackage{xcolor}
\usepackage{subfig}
\usepackage{authblk}
\usepackage[a4paper, margin=2cm]{geometry}

\makeatletter  

\numberwithin{equation}{section}
\newtheorem{theorem}{Theorem}
\newtheorem{definition}{Definition}
\newtheorem{proposition}{Proposition}
\newtheorem{lemma}{Lemma}
\newtheorem{corollary}{Corollary}
\newtheorem{remark}{Remark}
\newtheorem{claim}{Claim}

\newcommand{\mc}[1]{\mathcal{#1}}

\newcommand{\diag}{\textrm{diag}}
\newcommand{\sgn}{\textrm{sgn}}

\begin{document}

\title{Signed graphs in data sciences via communicability geometry}
\author[1]{Fernando Diaz-Diaz}
\author[1]{Ernesto Estrada}
\affil[1]{Institute of Cross-Disciplinary Physics and Complex Systems, IFISC (UIB-CSIC), 07122 Palma de Mallorca, Spain}
\date{}
\maketitle
\begin{abstract}
Signed graphs are an emergent way of representing data in a variety
of contexts where antagonistic interactions exist. These include data
from biological, ecological, and social systems. Here we propose the concept
of communicability for signed graphs {and explore in depth its mathematical properties. We also prove that the communicability induces a hyperspherical geometric embedding of the signed network, and derive communicability-based metrics that satisfy the axioms of a distance even in the presence of negative edges.} We then apply these metrics to solve several
problems in the data analysis of signed graphs within a unified framework. These
include the partitioning of signed graphs, dimensionality reduction,
finding hierarchies of alliances in signed networks, {and quantifying} the degree of polarization between the existing
factions in social systems represented by these types of graphs.
\end{abstract}
\maketitle

\section{Introduction}
Data emerging from a variety of contexts can be
represented in the form of signed graphs \cite{Zaslavsky1982a}.
A signed graph $\Sigma=(G,\sigma)$ consists of an underlying graph
$G=(V,E)$ and a signature function $\sigma:E\rightarrow\left\{ 1,-1\right\} $. The vertices $v_{i}\in V$ of $G$ represent entities in complex systems, {from molecules or species in biology to people in social systems}. 
The edges, representing the interactions between these entities, are
signed according to the nature of these interactions --with positive
signs reserved for activation, cooperation, friendship, alliances,
etc., and negative signs used to represent inhibition, competition,
foes, enmity, etc. For example, in a transcription network,
positive/negative edges represent the activation/inhibition activity
of a transcription factor on a target gene, while in metabolic and signaling networks, a signed
edge represents the contribution (accelerating or slowing
down) of a molecular species to a kinetic reaction \cite{iacono2010monotonicity}. In ecological
systems \cite{saiz2017evidence}, positive links model mutualistic
and facilitation interactions, while negative links
represent competitive and parasitic interactions. The interactions between the entities
of a complex system are, on many occasions, extracted from the real world
via correlations of various kinds. This is the case, for instance,
of financial networks \cite{macmahon2015community}
where correlations between stocks are frequently used to represent
the possible causal connections between the entities. Social systems are, by far, the richest sources of signed interactions
between entities forming a network. They include the cases of friends/foes
and collaborators/defectors in face-to-face and online social networks
\cite{leskovecSignedNetworksSocial2010}, the study of elections and voting systems in different scenarios \cite{arinik2020multiple},
and the evolution of international relations between countries during
regional/global conflicts \cite{diaz-diaz2023}.

The study of signed graphs is an active area of investigation in mathematics
\cite{Zaslavsky1982a}, with roots in the works of Harary in the 1950s
\cite{harary1953, cartwright1956structural}.
There have been many recent advances in the understanding of structural
properties of signed networks, as well as the dynamical processes taking
place on them (see e.g. \cite{tangSurveySignedNetwork2016}).
However, seeing these graphs as a way of representing complex data immediately
triggers a series of research questions. For instance, (i) how can
a signed graph be partitioned to detect the main factions existing
in the data; (ii) considering the signed graph as an $n$-dimensional
space of data, how can we reduce its dimensionality?; (iii) is there
a hierarchy of alliances in a signed graph?; (iv) how much polarization do
the existing factions have?; (v) can we predict
missing signed edges from the existing data? (see
the section ``Related work'' for details). One way to
collectively tackle these problems is by using machine learning techniques
resting on pairwise distances between data points. A distance metric
defines the geometry of an underlying space \cite{devriendt2022effective}, such as Euclidean, spherical,
or hyperbolic ones, which are at the core of many machine learning
algorithms \cite{tabaghi2020geometry}. The definition of such metrics
for signed graphs is the main goal of this work.

In this work, we define the concept of communicability geometry for
graphs with signs. Particularly, we define the communicability distance
and angle for signed graphs, which are based on the matrix exponential
of the signed adjacency matrix. We prove that these metrics are Euclidean
and spherical, showing several interesting properties for their use
in data analysis based on signed graphs. We apply the communicability
metrics to solve in a unified way the problems of signed graph partition,
dimensionality reduction, finding hierarchies of alliances in signed
networks as well as the quantification of the degree of polarization
between the existing factions in systems represented by this type
of graph. Although we do not tackle the problem of predicting missing
signed edges in signed graphs, the methods developed here can be applied
to solve this task based on similarities between pairs of vertices.
The solutions to all these problems are illustrated through real-world examples,
including the identification of conflicting tribal factions in Papua New Guinea, the analysis of World Wars I and II, the study of polarization within the European Parliament, {and the study of the gene regulatory network of \textit{E. coli}}.

\section{Definitions}  \label{sec:definitions}

Throughout this work, we will use the terms graph and network interchangeably.
We denote $\#|V|=:n$ and $\#|E|=:m$. We consider
only undirected connected signed graphs. Therefore, the adjacency
matrix $A\left(\Sigma\right)=A$ is a symmetric, square matrix whose entries
are $A_{ij}=\pm1$ if the nodes $i$ and $j$ are connected by a positive or negative edge, respectively; otherwise, $A_{ij}=0$.

We define the spectrum of $\Sigma$ as the set of eigenvalues of its adjacency
matrix: $Sp(\Sigma)=Sp(A)=\{\lambda_{1},...,\lambda_{n}\}$, where the eigenvalues
are in non-increasing order: $\lambda_{1}\geq\lambda_{2}\geq...\geq\lambda_{n}$.
The normalized eigenvector associated with the eigenvalue $\lambda_{i}$
is $\psi_{i}=(\psi_{i}(1),...,\psi_{i}(n))^{T}$, and $U=[\psi_{1}\ \psi_{2}\ \ ...\ \psi_{n}]$
denotes the orthogonal matrix whose columns are the eigenvectors of
$A$. To eliminate ambiguity regarding the direction of each eigenvector,
we consistently select $\psi_{i}$ so that its first component is
non-negative, i.e., $\psi_{i}(1)\geq0$ for every $i$. Additionally,
throughout all this work, $|A|$ denotes the entrywise absolute value
of $A$. $|A|$ defines an unsigned graph $G$,
called the \textit{underlying graph}, with the same sets of nodes
and edges as $\Sigma$ but where all the edges are positive.
A \textit{partition} $\mathcal{P}$ of the node set $V$ is a family
of subsets $V_{1},...,V_{n}$ that are pairwise disjoint and that
cover $V$. That is, $V_{i}\cap V_{j}=\emptyset$ for all $i,j$,
and $\bigcup_{j}V_{j}=V$. In the context of this work, we also
use the term \textit{factions} to name the subsets $V_{1},...,V_{n}$. When the set of nodes is partitioned into two factions,
$\mathcal{P}=\{V_{1},V_{2}\}$, we define an indicator function 
$I:V\to\{\pm1\}$ such that $I (v)=1 \textrm{ if }  v\in V_{1}$ and $I(v)=-1 \textrm{ if } v\in V_{2}$.

A \textit{walk} {$W$}
is a sequence of (not necessarily different) nodes $(v_{1},...,v_{n})$
such that consecutive nodes are connected by an edge: $(v_{i},v_{i+1})\in E$.
If the initial and final nodes coincide in a walk, $v_{1}=v_{n}$,
the walk is said to be \textit{closed}. A walk in which every node
and edge appears only once is called a \textit{path}. A closed path
is known as a \textit{cycle}. The \textit{sign} of the walk/path/cycle
is defined as the product of the signs of their edges: {$\sigma(W)=\prod_{e\in W} \sigma(e)$}. Consequently,
a \textit{positive} walk/path/cycle contains an even number of negative
edges, whereas a \textit{negative} one contains an odd number of negative
edges.

\begin{figure}[h]
\centering \includegraphics[width=0.8\textwidth]{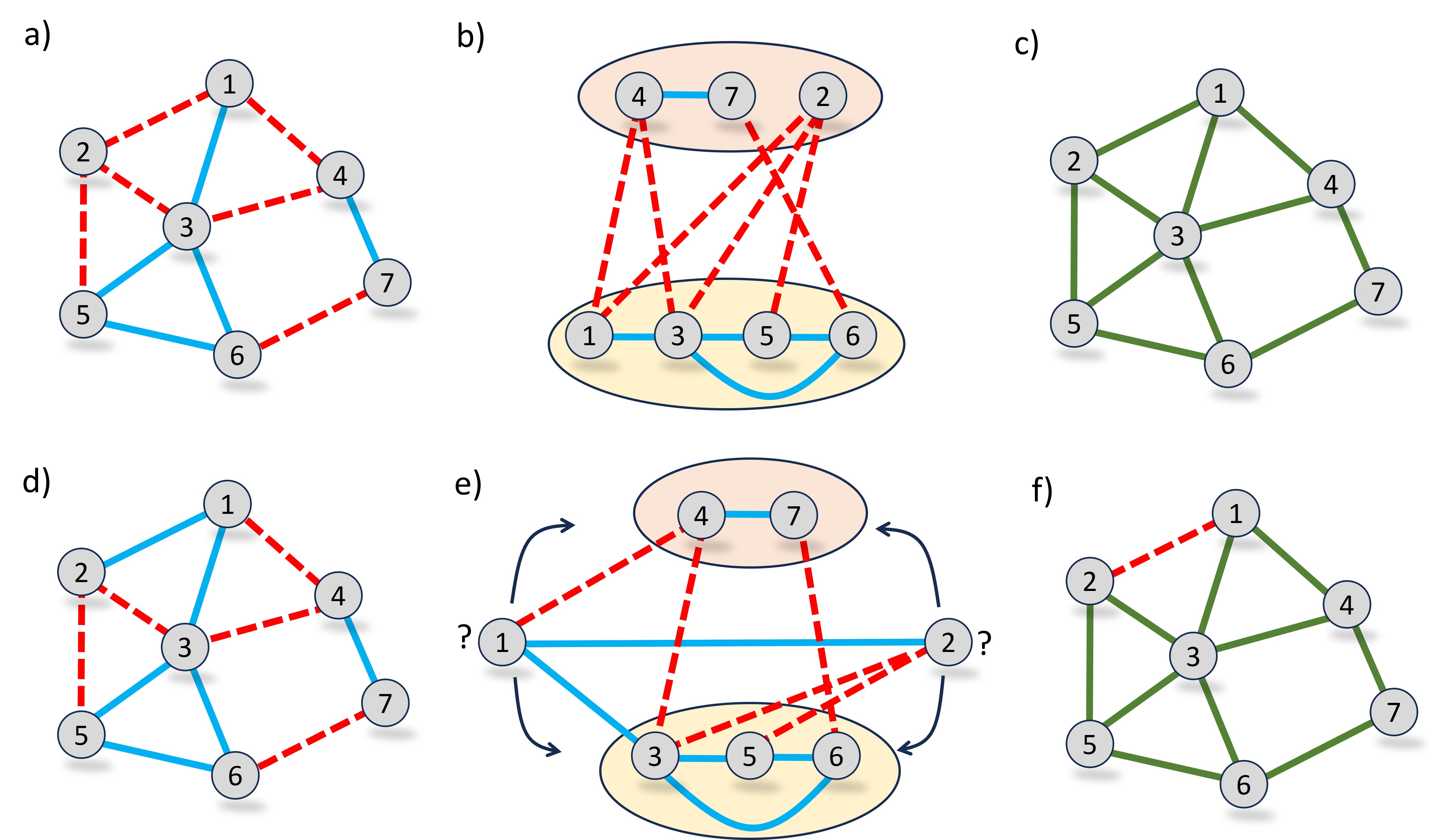} 
\caption{
Example of a balanced graph (a), its partition according to Theorem
{\ref{thm:Harary}} (b), and the result of applying the switching transformation described
in the text (c). Example of an unbalanced graph {(i.e. a graph containing cycles with negative parity)} (d), an illustration
of some problems emerging when trying to apply Theorem {\ref{thm:Harary}} to this
graph (e) as well as the graph resulting from applying the switching
transformation mentioned before (f). {Blue and green color denotes positive edges, while red dashed lines correspond to negative edges. }}
\label{fig:balance_examples}
\end{figure}

A signed graph is said to be \textit{balanced} \cite{harary1953} if every cycle within
it is positive or if it does not contain any cycle at all; or equivalently, if every closed
walk within it is positive (see Fig. \ref{fig:balance_examples}(a)).  The structural balance theorem proved by Harary \cite{harary1953} connects the concept of balance and that of a partition in a signed graph:

\begin{theorem}[Harary]  \label{thm:Harary}
Let $\Sigma$ be a signed graph. Then, $\Sigma$ is balanced
if and only if its node set admits a balanced bipartition; i.e., a
partition into \textit{balanced factions} $V=V_{1}\cup V_{2}$ such that every edge connecting nodes
of the same faction is positive, while every edge connecting nodes
of different factions is negative: 
$I\left(v_{i}\right)=I\left(v_{j}\right)$ if and only if $A_{ij}>0$
and $I\left(v_{i}\right)\neq I\left(v_{j}\right)$ if and only if
$A_{ij}<0$.
\end{theorem}

This theorem provides a second definition of balance, focusing on partitions instead of cycles. Indeed, in Fig. \ref{fig:balance_examples} we see that a balanced graph admits a balanced partition into two factions (panel (b)), whereas an unbalanced graph does not admit such a partition; i.e., we cannot
decide in which faction to place some of the nodes of unbalanced graphs (panel (e)).

A third way of defining a balanced graph is based
on the so-called \textit{switching transformations}. Given a graph
$\Sigma$ with adjacency matrix $A$, a graph switching is a
function $\mathcal{S}:\Sigma\to\Sigma'$ defined by $\mathcal{S}(A(\Sigma))=DA(\Sigma)D^{-1}$.
The matrix $D$ is a diagonal matrix whose diagonal elements are $+1$
or $-1$, often called a \textit{switching matrix}. Switching
matrices are symmetric and orthogonal; i.e., $D=D^{-1}=D^{T}$. If
two graphs $\Sigma$ and $\Sigma'$ can be related by a
switching transformation; i.e., $A(\Sigma')=DA(\Sigma)D^{-1}$,
we say that these graphs are\textit{ switching-equivalent}. Switching-equivalent
graphs have the same spectrum, because adjacency matrices related
by a switching transformation are similar. The following proposition
establishes a connection between structural balance, switching transformations,
and non-negative matrices: 
\begin{proposition}
\cite{hou2003} \label{switching_proposition} Let $\Sigma$ be a signed
graph with adjacency matrix $A$. $\Sigma$ is balanced if and
only if there exists a switching matrix $D$ such that $DAD^{-1}=|A|$.
Moreover, if such $D$ exists, it has the following form: $D=diag(I(v_{1}),...,I(v_{n}))$,
where $I(v_{j})$ is the indicator function defined before. 
\end{proposition}

In other words, the switching transformation maps a balanced graph into an
unsigned graph, see panel (c) of Fig. \ref{fig:balance_examples}.
On the other hand, if we apply the {same} switching matrix $D$ to the unbalanced
graph of panel (d), a negative edge appears between nodes 1 and 2
(see panel (f) in that Figure). {What's more, any choice of $D$ will result in a graph with at least one negative edge}. The appearance of this negative edge
is a direct consequence of the fact that edge (1,2) is 'frustrated',
i.e., we cannot decide a priori in which faction to place its endpoints in the partition given in panel (e).

\begin{remark}
A consequence of {Proposition \ref{switching_proposition}} is that all complete graphs
which are balanced are isospectral {(i.e. graphs whose adjacency matrices have identical spectra)}, so there are $n/2$ ($n$
even) or $\left(n-1\right)/2$ ($n$ odd) isospectral balanced complete graphs.
We illustrate in Fig. \ref{complete graphs} the three existing
ones with 6 nodes.

\begin{figure}[h]
\centering{}
\includegraphics[width=0.75\textwidth]{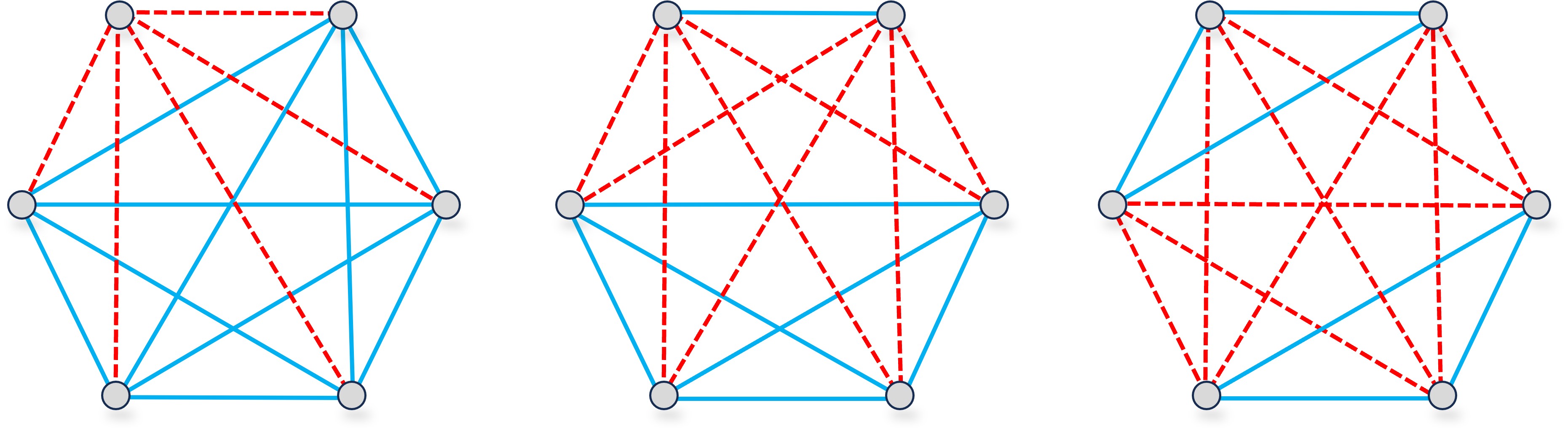}

\caption{Illustration of the three existing balanced complete graphs with 6
vertices, which are isospectral to one another.}

\label{complete graphs}
\end{figure}
\end{remark}

\section{Related works}  \label{sec:related_works}

Let us start by summarizing the efforts in the literature to consider distances
in signed graphs. Spiro \cite{spiroWienerIndexSigned2022} defined the length
of a path $P$ in $\Sigma$ as ${\tilde \sigma}(P)=\sum_{e\in {P}}\sigma\left(e\right)$,
where $\sigma:E\rightarrow\left\{ 1,-1\right\} $.
The shortest path ``signed distance'' between $u,v\in V$ was defined
as $d_{G,\sigma}\left(u,v\right):=\min_{P}\left|{\tilde \sigma}(P)\right|,$
where the minimum ranges over all $uv$-paths. Although this
``signed distance'' may be useful in some graph theoretic problems,
it is not suitable as a pairwise dissimilarity in signed graphs. For instance,
if $P=\left\{ \left(u,v\right),\left(v,w\right)\right\} $ where $\sigma\left(\left(u,v\right)\right)\neq\sigma\left(\left(v,w\right)\right),$
then $d_{G,\sigma}\left(u,w\right)=0$, in spite of the fact that
$u\neq w$. Thus, the ``signed distance'' is not even a distance.
More recently, Hameed et al. \cite{hameedSignedDistanceSigned2021}
defined two types of shortest-path distances in signed graphs by introducing
the following two auxiliary signs: (S1) $\sigma_{max}(u,v)=-1$ if
all shortest $uv$-paths are negative, and $+1$ otherwise; (S2) $\sigma_{min}(u,v)=+1$
if all shortest $uv$-paths are positive, and $-1$ otherwise. The
distances are then: 
\begin{enumerate}
\item $d_{max}(u,v)=\sigma_{max}(u,v)d(u,v)=\max\left\{ \sigma(P(u,v)):P(u,v)\in\mathbb{P}(u,v)\right\} d(u,v)$;
\item $d_{min}(u,v)=\sigma_{min}(u,v)d(u,v)=\min\left\{  \sigma(P(u,v)):P(u,v)\in\mathbb{P}(u,v)\right\} d(u,v)$,
\end{enumerate}
where $\mathbb{P}\left(u,v\right)$
is the collection of all shortest paths $P(u,v)$, and $d(u,v)$ is
the length of the shortest path in the underlying unsigned graph $G$.
{Clearly, this distance is not well-defined (in the sense that $d_{max}=d_{min}$) unless the graph is balanced, which is too restrictive for real datasets. Moreover, it is not positive and can violate the triangle inequality, making it a poor choice for a distance. These two examples show the difficulties of proposing a well-defined distance that incorporates the information contained in the negative edges while fulfilling the distance axioms. }

On the other hand, there have also been efforts in the literature to approach the problem
of signed graph partition from a computational perspective. The problem
is sometimes equated to the ``community detection'' one, typically
applied to unsigned networks, and consisting of the detection of densely
connected subsets of vertices in a graph \cite{fortunato2016community}. The signed version involves detecting subsets of nodes where positive edges connect nodes within each subset, while negative edges appear between subsets. The methods used include spectral techniques based on the adjacency matrix 
or the signed Laplacian \cite{kunegisSpectralAnalysisSigned2010,cucuringuSPONGEGeneralizedEigenproblem2019}, 
approaches that maximize modularity \cite{gomezAnalysisCommunityStructure2009}, and inferential methods based on the SBM \cite{jiang2015stochastic,peixoto2018nonparametric}, among others \cite{tangSurveySignedNetwork2016}.
The problem of dimensionality reduction in signed
graphs has also been addressed using embeddings based on the signed
Laplacian \cite{kunegisSpectralAnalysisSigned2010} or Hamiltonian functions \cite{babul2024sheep}.



\begin{figure}[h]
\begin{centering}
\includegraphics[width=0.98\textwidth]{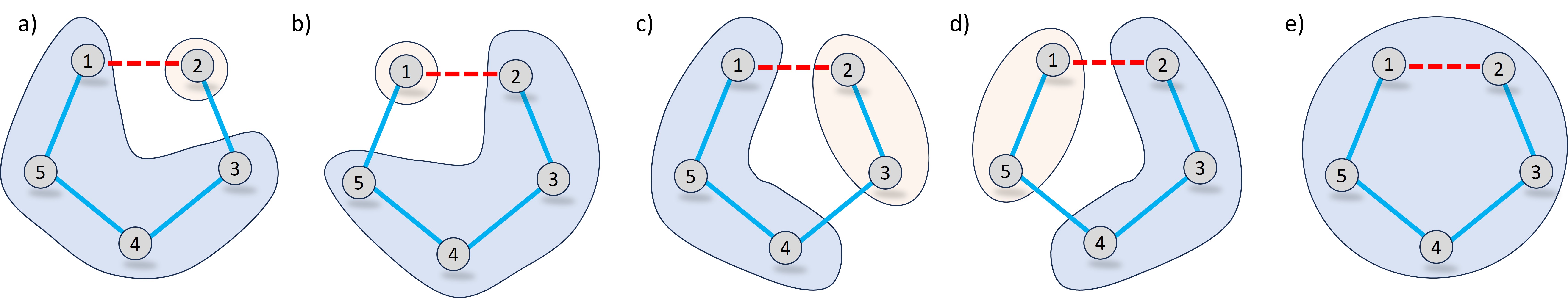}
\par\end{centering}
\caption{Example of five different signed partitions having the same ``frustration''
on the pentagon with one negative edge.}

\label{Example_1}
\end{figure}

If we closely examine the methods used for the analysis of
signed graphs, we find a series of drawbacks common to several
of them. To illustrate this, let us consider the
signed pentagon in which only one edge is negative. This example is
problematic for any method based on an objective function, such as
the frustration index. For instance, in Fig. \ref{Example_1} we give five
partitions of this graph that minimize the number of frustrated
edges. No other $k$-partitions minimize this objective
function. Thus, this problem is undecidable on the basis of frustration-based
minimization algorithms. Note that a frustration index based solely on triangles, as commonly defined, does not resolve this issue since the graph under study has no triangles.

Other approaches rely on either the leading eigenpair of the signed adjacency
matrix or the eigenpair corresponding to the smallest eigenvalue of
the signed Laplacian \cite{chiangScalableClusteringSigned2012}.
These methods are inspired by those that use such eigenvectors in
unsigned graphs \cite{vonluxburgTutorialSpectralClustering2007}.
While these eigenpairs are unique in unsigned graphs, this is not necessarily the case in signed graphs. Specifically, in a signed graph, the largest eigenvalue of $A$ may have multiplicity greater than one, which immediately introduces challenges for these methods. In the example
provided in Fig. \ref{Example_1} the eigenvalues are: $Sp\left(A\right)=\left\{ \varphi,\varphi,1-\varphi,1-\varphi,-2\right\} $,
where $\varphi=\cfrac{1+\sqrt{5}}{2}$ is the golden ratio. Therefore,
the largest eigenvalue is degenerate, i.e., it has a multiplicity two.
The two eigenvectors associated with the eigenvalue $\lambda_{1}$
of multiplicity two are: $\psi_{1}=\left[-1,\varphi,\varphi,1,0\right]^{T}$
and $\psi_{2}=\left[\varphi,-\varphi,-1\text{,0,1}\right]^{T}$, which correspond to the partitions (b) and (c) in Fig. \ref{Example_1}. The vector
$\psi_{1}+\varphi\psi_{2}=\left[\varphi,-1,0,1,\varphi\right]^{T}$,
which is an eigenvector of $A$ corresponding to $\lambda_{1}$
induces the partition (a) in Fig. {\ref{Example_1}}. Similarly, $-\psi_{1}=\left[1,-\varphi,-\varphi,-1,0\right]^{T}$
and $\varphi\psi_{1}+\psi_{2}=\left[0,1,\varphi,\varphi,1\right]^{T}$
induce the partitions (d) and (e), respectively. Therefore, the leading eigenvector of the signed
adjacency matrix does not induce a good partition when the associated eigenvalue has an algebraic multiplicity greater than one. Furthermore, for this graph, where the Laplacian is $L=2I-A,$ with $I$
being the identity matrix, there is a trivial relationship between the eigenvalues and eigenvectors of the Laplacian and those of the adjacency matrix. As a result, methods based on the smallest Laplacian eigenpair also fail to provide a solution for partitioning this graph.

Methods based on modularity--a measure of the quality of a partition
very popular in network analysis--also have several deficiencies, even
for unsigned graphs. There are other measures of the quality
of a partition, e.g., Silhouette, Davies-Bouldin, Calinski-Harabasz, etc. \cite{ pereda2019visualization}, and no systematic comparison between them and modularity
for (un)signed networks exists. Nonetheless, modularity has been generalized
to the case of signed networks
\cite{gomezAnalysisCommunityStructure2009}. This generalization
can identify the partition (a) in Fig. \ref{Example_1}
as the best one. However, it inherits some of the problems present in the unsigned case,
like the resolution
limit \cite{fortunatoResolutionLimitCommunity2007}. For instance,
let us build a signed graph formed by $k$ cliques labeled consecutively
by $1,2,\ldots,k$, each of size $r$. Two cliques $K_{i}$ and $K_{j}$
are connected by two edges, one positive and the other negative, if
and only if $j=i+1$ or $i=1$ and $j=k$. The obvious partition is
to consider every clique as one cluster. However, for sufficiently large $k$, the modularity approach identifies
the partition which groups pairs of cliques as the best one. This is
the signed version of the previously studied resolution problem in unsigned
graphs \cite{fortunatoResolutionLimitCommunity2007}.

\section{Signed walks, factions, and communicability}

In a social system, where the nodes of a network represent social
entities and edges represent their relations, a positive edge signifies an agreement, i.e. an affirmation of a node's state by its positively linked partner. Similarly, a negative edge signifies a conflict, i.e. a negation of the node's state by its partner. Let a double negation be considered as an
affirmation. Then, in a balanced signed graph, every walk starting at any
node $i$ in the balanced faction $V_{k}$ and ending at any node $j$ in the same faction represents an affirmation of the states
of the two nodes. Let us take for instance the pair of nodes 2 and
4 in graph (a) of Fig. \ref{fig:balance_examples}. The
state of node 4 is negated by node 1, while node 2 negates
the state of 1, which implies that nodes 2 and 4 affirm their states
to each other. On the other hand, the opposite behavior appears when the nodes belong to different (balanced) factions. Let us, for instance, consider the nodes 1
and 7 in the graph (a) of Fig. \ref{fig:balance_examples}.
In this case, the path $1\to4\to7$ results in a net negation between nodes 1 and 7, since node 4 negates the state of 1 and node 7, by affirming node 4's state, maintains that negation. In short, paths (or walks) negating the states of their endpoints are negative, while those mutually affirming the states of the two endpoints are positive.

In light of the previous discussion, we can say that \textit{two nodes $i$ and $j$
are effective allies if the number of positive walks connecting them
exceeds the number of negative walks; otherwise, the nodes are effective
enemies}. Moreover, the contribution of a walk should also depend on its length: a short
    (positive or negative) walk will have significantly more influence than a longer one. Mathematically,
we can express these ideas as follows. Let $\mu_{k}^{+}\left(i,j\right)$
be the number of positive walks of length $k$ connecting the nodes
$i$ and $j$, and let $\mu_{k}^{-}\left(i,j\right)$ be the number of
negative walks of length $k$ connecting the same nodes. Let

\begin{equation}
\Gamma_{ij}=\sum_{k=0}^{\infty}c_{k}\left(\mu_{k}^{+}\left(i,j\right)-\mu_{k}^{-}\left(i,j\right)\right),
\end{equation}
where $c_{k}$ is a penalization factor that decreases as $k$ increases.
Then, if $\Gamma_{ij}>0$ the nodes $i$ and $j$ are effective
allies, and if $\Gamma_{ij}<0$ they are effective enemies, and we have the following.
\begin{lemma} \label{switching_lemma}
Let $i,j\in V$ in a signed graph $\Sigma =(G,\sigma)$
and let \textup{$\mu_{k}^{+}\left(i,j\right)$ and }$\mu_{k}^{-}\left(i,j\right)$
be the number of positive and negative walks of length $k$ between
$i$ and $j$. Then,

\begin{equation}
\mu_{k}^{+}\left(i,j\right)-\mu_{k}^{-}\left(i,j\right)=\left(A^{k}\right)_{ij}.
\end{equation}
\end{lemma}

\begin{proof}{
    We proceed by induction on $k$. For $k=1$, there are only three possibilities: $i$ and $j$ are either connected by a positive edge, by a negative edge, or disconnected. 
    It is trivial to show that, in all three cases, it is true that  $\mu^+_1(i,j)-\mu^-_1(i,j)=A_{ij}$. \\
    Now, let us assume that $\mu^+_k(i,j)-\mu^-_k(i,j)=(A^k)_{ij}$ and prove it for $k+1$. The induction hypothesis on $k$ implies that
    \begin{align}
        \begin{cases}
           \mu^+_k(i,j)+\mu^-_k(i,j)= (|A|^k)_{ij}; \\ \mu^+_k(i,j)-\mu^-_k(i,j)= (A^k)_{ij}.
        \end{cases}
    \end{align}
   Solving the system, we find that 
   $ \mu^+_k(i,j) = \frac1{2}[(|A|^k)_{ij}+(A^k)_{ij}], \ \mu^-_k(i,j) = \frac1{2}[(|A|^k)_{ij}-(A^k)_{ij}] $. 
    
   We need to express the number of positive and negative walks of lengths $k+1$ in terms of the number of walks of length $k$. To do so, 
   we observe that each positive walk of length 
$k+1$ from $i$ to $j$ is composed of a walk of length $k$ from $i$ to some node $l$ and an edge $(l,k)$, with the $k$-walk and the last edge sharing signs. Similarly, a negative walk of length 
$k+1$ from $i$ to $j$ is composed of a walk of length $k$ from $i$ to some node $l$ and an edge $(l,k)$, with the $k$-walk and the last edge having  opposite signs. Thus, the total number of positive walks connecting $i$ and $j$ can be found by summing over the neighbors of $j$, distinguishing the positive neighbors, $\mc N ^+(j)$, from the negative ones, $\mc N^-(j)$:
     \begin{align}
        \mu^+_{k+1}(i,j) = \sum_{l\in\mc N^+(j)} \mu^+_{k}(i,l) +  \sum_{l\in\mc N^-(j)} \mu^-_{k}(i,l).
    \end{align}
   The condition $l\in\mc N^{\pm}(j)$ is equivalent to multiplying by $\mu_1^{\pm}(l,j)$; therefore:  
    \begin{align}
        \mu^+_{k+1}(i,j) = \sum_{l=1}^N \mu^+_{k}(i,l)\mu^+_{1}(l,j) +  \sum_{l=1}^N \mu^-_{k}(i,l)\mu^-_{1}(l,j).
    \end{align}
Now we can apply the expressions for $\mu^\pm_k(i,j)$ to get:
    \begin{align}
        \mu^+_{k+1}(i,j) &= 
         \sum_{l}  \left(\frac{(|A|^k)_{il}+(A^k)_{il}}{2}  \right) \left(\frac{|A|_{lj} + A_{lj}}{2} \right) + \sum_l \left( \frac{(|A|^k)_{il} -  (A^k)_{il}}{2} \right) \left(\frac{|A|_{lj} - A_{lj}}{2} \right) = \nonumber \\
        &= \frac1{2} \sum_l [(|A|^k)_{il}|A|_{lj} + (A^k)_{il}A_{lj}] 
         = \frac1{2} [(|A|^{k+1})_{ij} + (A^{k+1})_{ij}].  \nonumber
    \end{align}
In the same manner, $\mu^-_{k+1}(i,j)$ can be written as:
\begin{align}
     \mu^-_{k+1}(i,j) &= \sum_{l=1}^N \mu^+_{k}(i,l)\mu^-_{1}(l,j) +  \sum_{l=1}^N \mu^-_{k}(i,l)\mu^+_{1}(l,j) \nonumber \\ &= \frac1{4} \sum_{l}  \left(\frac{(|A|^k)_{il}+(A^k)_{il}}{2} \right) \left(\frac{|A|_{lj} - A_{lj}}{2} \right) + \sum_l\left(\frac{(|A|^k)_{il}-(A^k)_{il}}{2} \right) \left(\frac{|A|_{lj} + A_{lj}}{2} \right) = \nonumber \\
     &= \frac1{2} \sum_l [(|A|^k)_{il}|A|_{lj} - (A^k)_{il}A_{lj}] =  \frac1{2}[(|A|^{k+1})_{ij} - (A^{k+1})_{ij}].  \nonumber
\end{align}
    Finally, subtracting:
    \begin{align}
        \mu^+_{k+1}(i,j) - \mu^-_{k+1}(i,j) &= \frac1{2} [ (|A|^{k+1})_{ij}+ (A^{k+1})_{ij} - (|A|^{k+1})_{ij} + (A^{k+1})_{ij} ] = (A^{k+1})_{ij}, \nonumber
    \end{align}
    which is the induction hypothesis for $k+1$. Thus, by the induction principle, our claim is true for all $k$.
    }
\end{proof}

Let us now consider the following.
\begin{definition} \label{c-communicbaility}
Let $\Sigma$ be a signed graph with adjacency matrix $A$. The \textit{c-signed
communicability} $\Gamma_{ij}(A)$ function between nodes $i$
and $j$ is the $\left(i,j\right)$-entry of the matrix function $f\left(A\right)$
given by: 

\begin{equation}
\Gamma_{ij}(A)= \sum_{k=0}^{\infty}c_{k}\left(A^{k}\right){}_{ij}=\left(f\left(A\right)\right)_{ij}.
\end{equation}
\end{definition}

{Throughout this work,} we focus on the factorial penalization $c_{k}=\left(k!\right)^{-1}$ and
 henceforth refer to $\Gamma_{ij}$ as the signed communicability,
even though it corresponds to the exp-signed communicability. We note that this choice of penalization has been shown to be useful for defining indices that quantify both global \cite{estrada2014balance} and local \cite{diaz-diaz2023} balance, as well as for measuring polarization in signed graphs \cite{talagaPolarizationMultiscaleStructural2023}.
 Other penalization factor can nevertheless be selected to give more or less weight to
longer walks, resulting in other matrix
functions \cite{estrada2022many}.

The term $\Gamma_{ii}$ accounts for the signed self-communicability
of node $i$. Notice that in the argot of network theory, this
term represents a centrality measure of the given node counting its
participation in all signed subgraphs of the graph and giving more
weight to the shorter than to the longer ones. The following is a
difference between the self-communicability and the communicability
between two nodes in a signed graph.
\begin{lemma}
Let $\Sigma=(G,\sigma)$ be a connected signed
graph with adjacency matrix $A$. Let $\Gamma_{ij}(A)$ be the signed self-communicability between
the nodes $i$ and $j$. Then, $\Gamma_{ii}(A)>0$ for all $i\in V$, while
$\Gamma_{ij}(A)$ can take both positive and negative values.
\end{lemma}

\begin{proof}
{Since $A$ is symmetric, we can use its eigendecomposition: $A_{ij}=\sum_{v=1}^{n}{\lambda_{v}}\psi_{v}\left(i\right)\psi_{v}\left(j\right).$ Consequently, the exponential can be applied to the eigenvalues while leaving the eigenvectors unchanged, resulting in 
$\Gamma_{ij}(A)=\sum_{v=1}^{n}e^{\lambda_{v}}\psi_{v}\left(i\right)\psi_{v}\left(j\right).$ 
When $i=j$, this formula reduces to} $\Gamma_{ii}(A)=\sum_{v=1}^{n}e^{\lambda_{v}}\psi_{v}^{2}\left(i\right)>0,$
while when $i\neq j$ it can also take negative values depending on
the sign of $\psi_{v}\left(i\right)$ and $\psi_{v}\left(j\right)$.
\end{proof}

We now arrive at the following result that connects the signed and
unsigned communicability.
\begin{proposition}
\label{comm_bound} Let $\Sigma$ be a signed graph with adjacency
matrix $A$. Let $\Gamma_{ij}(A)$ be the communicability between
the nodes $i$ and $j$ in the signed graph, and let $\Gamma_{ij}(|A|)$ be
the communicability in the unsigned underlying graph of $\Sigma$.
Then, $|\Gamma_{ij}(A)|\leq\Gamma_{ij}(|A|)$. 
\end{proposition}
\begin{proof}

Let us write 
\begin{align}
|\Gamma_{ij}(A)|=|\left(e^{A}\right){}_{ij}|=\left|\sum_{k=0}^{\infty}\frac{1}{k!}(A^{k})_{ij}\right|\leq\sum_{k=0}^{\infty}\frac{1}{k!}|A^{k}|_{ij}.
\end{align}

On the other hand: 
\begin{align}
\left|A^{k}\right|_{ij}=\left|\sum_{l_{1},...,l_{k-1}}A_{il_{1}}A_{l_{1}l_{2}}...A_{l_{k-1}j}\right|\leq\sum_{l_{1},...,l_{k-1}}|A|_{pl_{1}}|A|_{l_{1}l_{2}}...|A|_{l_{k-1}j}=|A|_{ij}^{k}.\label{adjacency_inequality}
\end{align}
Combining the last two inequalities, we obtain $|\Gamma_{ij}(A)|\leq\sum_{k}\frac{1}{k!}|A|_{ij}^{k}=\Gamma_{ij}(|A|)$. 
\end{proof}

\subsection{Communicability of balanced graphs}

We start this subsection by proving the following result.
\begin{lemma}
\label{comm_balance_lemma} 
Let $\Sigma,\Sigma'$ be signed
graphs with adjacency matrices $A,A'$ respectively. Moreover, assume
that $\Sigma$ and $\Sigma'$ are switching equivalent;
i.e., there is a switching matrix $D$ such that $A'=DAD^{-1}$. Then,
the communicability matrices of $\Sigma$ and $\Sigma'$
are related by: 
\begin{align}\label{eq_comm_balance_theorem}
\Gamma(A')=D\Gamma(A)D^{-1}.
\end{align}
\end{lemma}
\begin{proof}
    Let us notice that $(A')^{2}=DAD^{-1}DAD^{-1}=DA^{2}D^{-1}$. Proceeding
    by induction, we find that $(A')^{k}=DA^{k}D^{-1}$. Therefore: 
    \begin{align}
    \Gamma(A')=e^{A'}=\sum_{k=0}^{\infty}\frac{(A')^{k}}{k!}=\sum_{k=0}^{\infty}\frac{DA^{k}D^{-1}}{k!}=D\sum_{k=0}^{\infty}\frac{A^{k}}{k!}D^{-1}=De^{A}D^{-1}=D\Gamma(A)D^{-1}.
    \end{align}
\end{proof}

\begin{theorem}
\label{comm_balance_theorem} Let $\Sigma$ be a signed graph
with adjacency matrix $A$. Then, $\Sigma$ is balanced if and
only if there exists a switching matrix $D$ such that: 
\begin{align}
\Gamma(A)=D\Gamma(|A|)D^{-1},\label{eq_comm_balance_theorem_2}
\end{align}
where $D$ is the switching matrix that maps $A$ to $|A|$.
\end{theorem}

\begin{proof}
{
We begin by proving the necessary condition. Assume that $\Sigma$
is balanced. Then, by Proposition \ref{switching_proposition}, there is
some matrix $D$ such that $A=D|A|D^{-1}$. Then, according to Lemma
\ref{comm_balance_lemma}, we have that $\exp(A)=D\exp(|A|)D^{-1}$. \\
 Let us now focus on the sufficient condition. For that let us assume
that $\Gamma(A)=D\Gamma(|A|)D^{-1}$, or equivalently, $\exp(|A|)=D^{-1}\exp(A)D$.
We then proceed to prove that $\Sigma$ is balanced. Let us use
the Taylor series expansion of the communicability, to get: 

\begin{equation}
\sum_{k=0}^{\infty}\frac{|A|^{k}}{k!}=\sum_{k=0}^{\infty}\frac{D^{-1}A^{k}D}{k!} \Rightarrow \sum_{k=0}^{\infty}\frac{1}{k!}(|A|^{k}-D^{-1}A^{k}D).
\end{equation}

Using that $D$ is diagonal and $D=D^{-1}$,
we can write each term within the parenthesis as $|A|_{ij}^{k}-D_{ii}D_{jj}(A^{k})_{ij}$.
From the proof
of Proposition \ref{comm_bound}, we get that $|A|_{ij}^{k}\geq|A^{k}|_{ij}\geq (A^{k})_{ij}$, so all terms within the sum are greater or equal to zero. Consequently,
the sum is zero if and only if $|A|^{k}=D^{-1}A^{k}D,$
for every $k$. In particular, taking $k=1$, we get $|A|=DAD^{-1},$ which according to Proposition \ref{switching_proposition} is a sufficient
condition for $\Sigma$ to be a balanced graph. Moreover, this
also proves that $D$ is indeed the balanced switching matrix of $A$,
completing the proof. }
\end{proof}

Let us now explore the many consequences of Theorem \ref{comm_balance_theorem}.
To begin with, let us go back to the inequality of Proposition \ref{comm_bound}:
$|\exp(A)|\leq\exp(|A|)$. Theorem \ref{comm_balance_theorem} allows
us to prove that the upper bound is reached when
the network is balanced: 
\begin{corollary}   \label{corollary_comm_balance}
 Let $\Sigma$ be a signed
graph with adjacency matrix $A$. $\Sigma$ is balanced if and
only if $|e^{A}|=e^{|A|}$. 
\end{corollary}

\begin{proof}
{
If $\Sigma$ is balanced, then Theorem \ref{comm_balance_theorem}
proves that $(e^{A})_{ij}=D_{ii}(e^{|A|})_{ij}D_{jj}$. Taking absolute
values on both sides and using that $|D_{ii}|=1$ for all $i$, we
get $|e^{A}|=e^{|A|}$.
On the other hand, if $\Sigma$ is unbalanced, then there are two nodes
$i$ and $j$ connected through both a positive and a negative path
\cite{harary1953}. Let us denote the lengths of the positive and
negative paths by $l_{+}$ and $l_{-}$, respectively. Recall that
$A^{k}$ counts the difference between positive and negative walks
of length $k$, while $|A|^{k}$ counts the number of walks, regardless
of sign. Moreover, Eq. (\ref{adjacency_inequality}) implies that
$-|A|_{ij}^{k}\leq A_{ij}^{k}\leq|A|_{ij}^{k}$ . The first inequality
is reached if and only if all paths between $i$ and $j$ are negative,
whereas the second one is reached if and only if all such paths are
positive. Therefore, the existence of a negative path of length $l_{-}$
between nodes $i$ and $j$ implies that $A_{ij}^{l_{-}}<|A|_{ij}^{l_{-}}$;
and the existence of the positive path of length $l_{+}$ implies
$A_{ij}^{l_{+}}>-|A|_{ij}^{l_{+}}$. This means that: 

\begin{equation}
-\sum_{k=0}^{\infty}\frac{|A|_{ij}^{k}}{k!}<\sum_{k=0}^{\infty}\frac{A_{ij}^{k}}{k!}<\sum_{k=0}^{\infty}\frac{|A|_{ij}^{k}}{k!},
\end{equation}
which implies that $-(e^{|A|})_{ij}<(e^{A})_{ij}<(e^{|A|})_{ij}.$
Thus, $|e^{A}|\neq e^{|A|}$, completing the proof. 
}
\end{proof}

\begin{corollary} \label{corollary_2}
Let $\Sigma$ be a signed complete graph. Then, $\Sigma$ is balanced if and
only if for any pair of different vertices $i$ and $j$, $\Gamma_{ij}(A)=\dfrac{e^{n}-1}{ne}$
or $\Gamma_{ij}(A)=\dfrac{1-e^{n}}{ne}.$ 
\end{corollary}

\begin{proof}
{

We begin proving the necessary condition. If $\Sigma$ is balanced,
it admits a balanced switching matrix $D$, and we can apply Eq. 
(\ref{eq_comm_balance_theorem_2}). Moreover, $i$ and $j$ belong to different
communities, so the elements $D_{ii}$ and $D_{jj}$ have opposite
sign. Thus, $\Gamma_{ij}(A)=-\Gamma_{ij}(|A|)=-\frac{e^{n}-1}{ne}$. \\
 We now prove the sufficient condition, which is equivalent to proving
that if $\Sigma$ is unbalanced, then $|\Gamma_{ij}(A)|\neq\dfrac{e^{n}-1}{ne}$.
If $\Sigma$ is unbalanced, then according to Corollary \ref{corollary_comm_balance}, there
is a pair $(i,j)$ such that $|\Gamma_{ij}(A)|<\Gamma_{ij}(|A|)$. Since $\Sigma$
is complete, it follows that $\Gamma_{ij}(|A|)=\frac{e^{n}-1}{ne}$. Thus,
$|\Gamma_{ij}(A)|\neq\frac{e^{n}-1}{ne}$. }
\end{proof}

\section{Communicability geometry}

As we have analyzed in Section \ref{sec:related_works}, the signed versions of the shortest path distance exhibit numerous problems, making them unsuitable as similarity measures for data analysis of signed graphs.
Thus, it is of great importance to find a distance function that adequately
captures the information contained in the negative edges while fulfilling
the distance axioms.

Let us start by considering a balanced partition of a balanced signed
graph $\Sigma$ in which $i\in V_{1}$ and $j\in V_{2}$.  Since $\Sigma$ is balanced, every closed walk starting (and ending)
at node $i$ (or at node $j$) is always positive.
Similarly, every walk starting at node $i$ and ending at node $j$
is negative. Consequently, we have $\mu_{k}^{-}\left(i,i\right)=\mu_{k}^{-}\left(j,j\right)=0$, as well as $\mu_{k}^{+}\left(i,j\right)=\mu_{k}^{+}\left(j,i\right)=0$. Furthermore, if $\Sigma$ is unbalanced but ``close to balance" (i.e. containing only a small number of negative cycles), these terms will be nonzero but close to zero for a suitably chosen partition of the node set into factions. Summing all these terms while accounting for their signs, and assuming that $\mu_{k}^{\pm}\left(i,j\right)=\mu_{k}^{\pm}\left(j,i\right)$
due to the the graph's undirected nature, we obtain the quantity 
${Q_k}= \left(\mu_{k}^{+}\left(i,i\right)-\mu_{k}^{-}\left(i,i\right)\right)+\left(\mu_{k}^{+}\left(j,j\right)-\mu_{k}^{-}\left(j,j\right)\right)-2\left(\mu_{k}^{+}\left(i,j\right)-\mu_{k}^{-}\left(j,i\right)\right)$, which is expected to be well-defined and positive. {Likewise, the sum $\sum_k \frac{Q_k}{k!} =\Gamma_{ii}+\Gamma_{jj}-2\Gamma_{ij}$ will also be positive}. This motivates the introduction of the following concept based on the communicability function (which,
as we have seen, accounts for a weighted difference between positive and
negative walks in $\Sigma$):
\begin{definition}
Let $\Gamma:=\exp\left(A\right)$ be the communicability matrix
of an undirected signed graph with adjacency matrix $A$. Then, let
us define the term $\xi_{ij}$ between nodes $i$ and $j$ as 
\begin{align}  \label{eq:comm_distance}
\xi_{ij}=\left(\Gamma_{ii}+\Gamma_{jj}-2\Gamma_{ij}\right)^{1/2}.
\end{align}
\end{definition}
We then have the following:
\begin{theorem}
\label{distance} $\xi_{ij}$ is a Euclidean distance between the
 nodes in the signed graph. 
\end{theorem}
\begin{proof}
{

We begin the proof by considering the spectral decomposition of $\Gamma_{ij}$:
$\Gamma_{ij}=\sum_{k}\psi_{k}(i)\psi_{k}(j)e^{\lambda_{k}}.$
Substituting this expression in the definition of $\xi_{ij}$ and
squaring both sides, we get $\xi_{ij}^{2}  =\sum_{k}e^{\lambda_{k}}\left(\psi_{k}(i)-\psi_{k}(j)\right)^{2}.$
We now define the vector $x_i$ as $x_{i}(k):=e^{\lambda_k/2}\psi_k(i)$. Then, we can write:
\begin{align}
    \xi_{ij}^2 = \sum_k \left(x_i(k)-x_j(k)\right)^2 = \Vert x_i - x_j\Vert^2,  \label{aux:position_vector_communicability}
\end{align}
where $\Vert \cdot \Vert$ is the norm induced by the usual inner product. Therefore,
we have shown that $\xi_{ij}$ is the Euclidean distance between the vectors $x_i$ and $x_j$.
}
\end{proof}

{We note that this metric is a natural extension of the communicability distance for unsigned networks \cite{estrada2012}, which has been extensively used in the literature (e.g. \cite{estrada2023network}).}
 As usual in the theory of Euclidean distance geometry, we define the
communicability Euclidean Distance Matrix (EDM) as the matrix whose
entries are the squares of the communicability distances:

\begin{equation}
M=s\mathbf{1}^{T}+\mathbf{1}s^{T}-2e^{A},
\end{equation}
where $\mathbf{1}$ is an all-ones column vector and $s:=\diag(e^A)$. We now state that
the communicability EDM is spherical \cite{estrada2014embedding}. An EDM is called spherical or circum-EDM if the set of points it represents lies on a hypersphere. 
\begin{lemma}
$M$ is a non-singular spherical EDM.
\end{lemma}
\begin{proof}
{
    The theorem is a special case of Theorem 5.1 of \cite{estrada2014embedding}, where the 
    authors prove that any EDM of the form $M=\vec 1^Ts +s^T\vec 1 - 2f(A)$, 
    with ${s}=\diag(f(A))$ and $f(A)$ a positive-definite matrix function, is a circum-Euclidean distance matrix.
    For the specific case of the signed communicability, it is clear that $f(A)=e^A$ is a valid positive definite matrix function even for signed adjacency matrices, proving the result.}
\end{proof}

\begin{figure}
\begin{centering}
\includegraphics[width=0.8\textwidth]{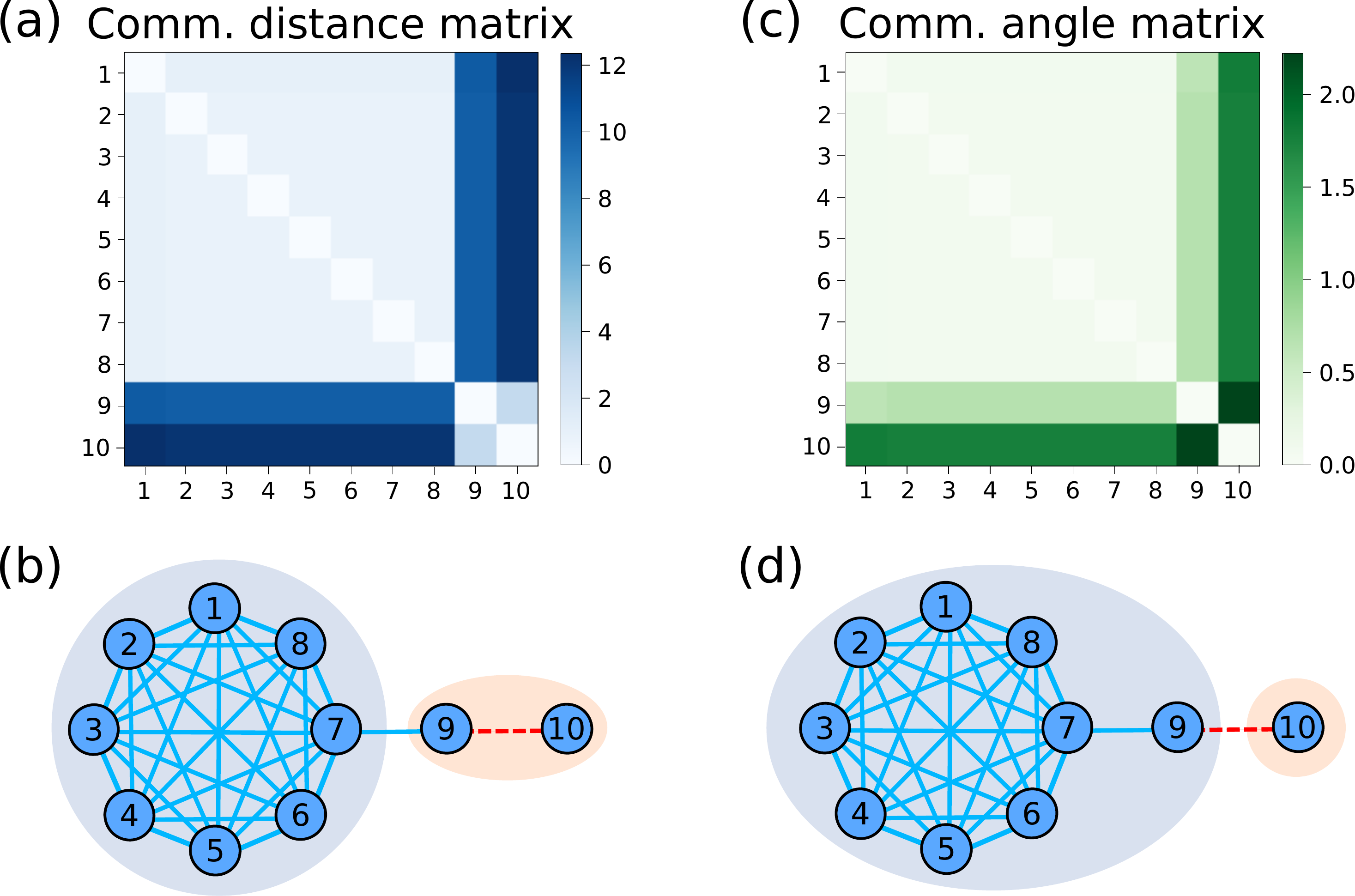}
\par\end{centering}
\caption{{Example illustrating a case where the algorithms based on the
communicability distance fail to identify the expected balanced partition.
    (a) Communicability distance matrix of the graph shown in (b). (b) Partition obtained by hierarchical clustering on the communicability distance matrix.
    Any clustering algorithm based on the communicability distance will fail to identify the 
    balanced partition, since the communicability distance between nodes 9 and 10 is smaller than between node 9 and any node in the clique.
    (c) Communicability angle matrix of the graph. (d) Partition obtained by hierarchical clustering on the communicability angle matrix. 
    Angle-based clustering correctly identifies the balanced partition.}
 }
\label{counterexample}
\end{figure}

A potential problem of the (signed) communicability distance is that
it is unbounded. {For instance, a direct calculation shows that the communicability distance between two nodes of a balanced complete graph belonging to opposite balanced factions is $\xi_{ij} = \sqrt{\frac{2(2e^n+n-2)}{ne}}$, which becomes arbitrarily large as $n\to\infty$. }
{Moreover,} the terms
$\Gamma_{ii}$ and $\Gamma_{jj}$ {introduce a bias} in the definition
of this distance. Let us provide an example. In the signed graph illustrated in Fig. \ref{counterexample}, we see that, in the balanced partition, node 10 belongs to one balanced faction, while the rest of the network forms the other one. However, based on the
signed communicability distance, nodes 9 and 10 are closer to each other than to any node within the clique {(Fig. \ref{counterexample}, panel (a))}. Therefore, any clustering algorithm based on the communicability distance will obtain
one cluster formed by nodes 1-8 and another formed by nodes 1-2 (panel (b)), which does not correspond to the balanced partition.

In the following, we will investigate the terms $\Gamma_{ii}$ and $\Gamma_{jj}$,
which are the main cause of the bias observed in the (signed) communicability
distance. We start by stating the following result.
\begin{proposition}
\label{position_vectors_proposition} Let $\{x_{i}\}$ be the generator set
formed by the communicability position vectors $x_{i}$, where $x_{i}(k)=\psi_{k}(i)e^{\lambda_{k}/2}$,
and $(\lambda_{k},\psi_{k})$ is the $k$th eigenpair of the adjacency
matrix $A$ of the graph $\Sigma$. Moreover, let $\hat{x}_{i}$
be the communicability position vectors of the underlying unsigned
graph and $I(i)$ is the indicator function defined in Section \ref{sec:definitions}.
Then, the following assertions hold: 
\begin{enumerate}
\item The squared norm of $x_{i}$ fulfills
$||x_{i}||^{2}=\Gamma_{ii}$. Moreover, if $\Sigma$ is balanced,
then $||x_{i}||=||\hat{x}_{i}||$. 
\item The inner product of $x_{i}$ and $x_{j}$ is the signed communicability:
$x_{i}\cdot x_{j}=\Gamma_{ij}$. 
\item $\{x_{i}\}$ is a basis of an $n$th-dimensional Euclidean space,
where $n$ is the number of nodes. 
\item There is only one $i$ for which $x_{i}$ contains only positive entries. 
\item If $\Sigma$ is balanced, then $x_{i}=I(i)\hat{x}_{i}$. 
\item If $\Sigma$ is balanced, then the sign of the first coordinate
of every $x_{i}$ gives the balanced partition of the graph. 
\end{enumerate}
\end{proposition}

\begin{proof}
{
The first two assertions come from the fact that 

\begin{align}
x_{i}\cdot x_{j}=\sum_{k}\psi_{k}(i)e^{\lambda_{k}/2}\psi_{k}(j)e^{\lambda_{k}/2}=\Gamma_{ij}.
\end{align}

In particular, applying this identity to the underlying graph, we
get $||\hat{x}_{i}||^{2}=\Gamma_{ii}(|A|)$. Finally, since $\Sigma$
is balanced, we apply Theorem \ref{comm_balance_theorem} to the diagonal
entries of $\Gamma$, obtaining $\Gamma_{ii}(A)=\Gamma_{ii}(|A|)$; i.e, $||x_{i}||=||\hat{x}_{i}||$. 

To prove statement 3, let us define $X$ as the matrix whose
columns are the vectors $x_{j}$. It is easy to see that $X=e^{\Lambda/2}U^{T}$:
\begin{align}
X_{ij}=x_{j}(i)=\psi_{i}(j)e^{\lambda_{i}/2}=\sum_{k}\delta_{ik}e^{\lambda_{k}/2}\psi_{k}(j)=\sum_{k}(e^{\Lambda/2})_{ik}U_{jk}=(e^{\Lambda/2}U^{T})_{ij}
\end{align}
Because neither $e^{\Lambda/2}$ nor $U$
are singular, the columns of $X$ form a set of linearly independent vectors. Since $\dim(X)=n$,
this set of vectors is a basis. 

We now prove the fourth statement. First we prove that at least one $x_{i}$ contains only positive
components. Since the eigenvectors $\psi_{j}$ of $A$ are defined
up to sign, one can choose the sign of each eigenvector so that
the $i$th component of every eigenvector is positive. Therefore,
$x_{i}=(\psi_{1}(i)e^{\lambda_{1}/2},\psi_{2}(i)e^{\lambda_{2}/2},...)^{T}$
only contains positive entries. To prove that there cannot be two
vectors $x_{i},x_{j}$ with only positive components, first observe
that any set of vectors of the form $(\psi_{1}(i),\psi_{2}(i),...)^{T}$
forms an orthonormal basis. Suppose $x_{i}$ contains only positive
entries. Then, $(\psi_{1}(i),\psi_{2}(i),...)^{T}$ also contains
only positive entries. Hence, to ensure orthogonality, every vector
of the form $(\psi_{1}(j),\psi_{2}(j),...)^{T}$ with $j\neq i$ must
contain at least one negative entry. The same entry will also be negative
in the vector $(\psi_{1}(i)e^{\lambda_{1}/2},\psi_{2}(i)e^{\lambda_{2}/2},...)^{T}$.
But this vector is $x_{j}$, so every $x_{j}$ with $j\neq i$ contains
a negative entry. 

To prove the fifth statement let us begin by relating the eigenvectors
of $A$ and $|A|$ when $\Sigma$ is balanced. Since $\Sigma$
is balanced, Proposition \ref{switching_proposition} relates $A$ and $|A|$
through $|A|=DAD$. Let $\hat{\psi}_{j}$ be an eigenvector of $|A|$.
Then, $\lambda_{j}\hat{\psi}_{j}=|A|\hat{\psi}_{j}=DAD\hat{\psi}_{j}\Rightarrow A(D\hat{\psi}_{j})=\lambda_{j}(D\hat{\psi}_{j})$;
i.e., $D\hat{\psi}_{j}$ is an eigenvector of $A$ with the same eigenvalue
as $|A|$. In other words, the eigenvectors $\psi_{j}$ and $\hat{\psi}_{j}$
are related by $\psi_{j}=D\hat{\psi}_{j}$, and their components fulfill
$\psi_{j}(i)=D_{ii}\hat{\psi}_{j}(i)$. 
 Now we substitute this last expression into the definition of $x_{i}$:
\begin{equation}
x_{i}(j)=\psi_{j}(i)e^{\lambda_{j}/2}=D_{ii}\hat{\psi}_{j}(i)e^{\lambda_{j}/2}=D_{ii}\hat{x}_{i}(j),
\end{equation}
where $\hat{x}_{i}$ is the communicability position vector of $|A|$.
Since this last equation is valid for all $j$ and $D_{ii}=I(i)$,
we conclude that $x_{i}=I(i)\hat{x}_{i}$. 

Finally, the sixth statement is a consequence of the fact that since
$|A|$ is non-negative, the Perron-Frobenius theorem asserts that
the principal eigenvector of $\hat{\psi}_{1}$ only contains positive
entries. Consequently, $\hat{x}_{i}(1)=\hat{\psi}_{1}(i)e^{\lambda_{1}/2}$
is positive. Now, because $\Sigma$ is balanced, we can use statement
5 to relate $x_{i}(1)$ and $\hat{x}_{i}(1)$ through $x_{i}(1)=I(i)\hat{x}_{i}(1)$,
and therefore $\sgn(x_{i}(1))=I(i)$. In other words, the sign of $x_{i}(1)$
tells us the partition subset to which node 1 belongs.  }
\end{proof}

\begin{figure}[h]
\centering
\includegraphics[width=0.8\textwidth]{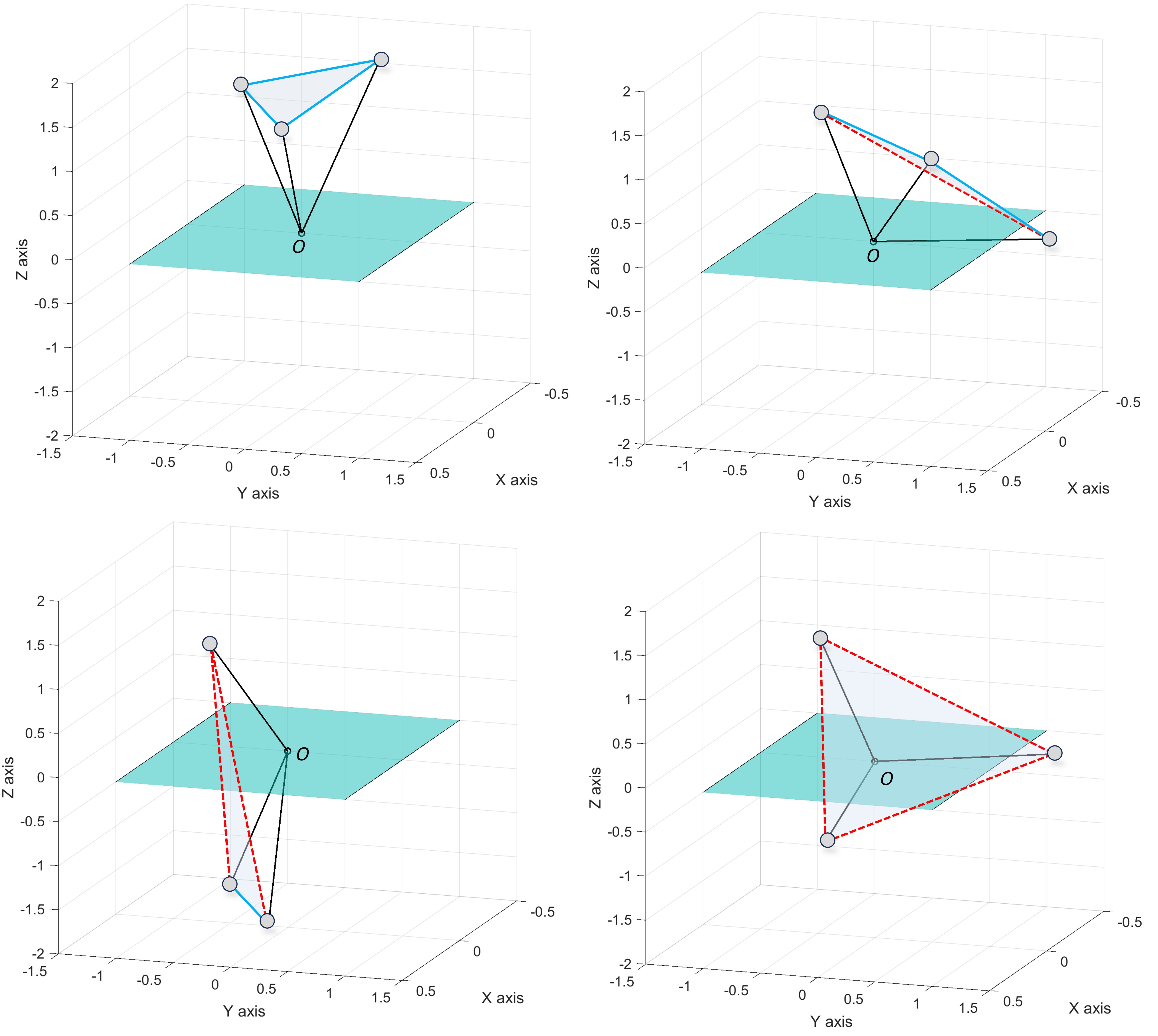}
\caption{Embedding of the four signed triangles in the induced three-dimensional communicability space.}
\label{triangles}
\end{figure}

One important consequence of Proposition \ref{position_vectors_proposition}
is that the position vectors contain important information to cluster
the nodes into factions, even in the unbalanced case. This is exemplified
in Fig. \ref{triangles}, where we plot the four signed triangles
in the induced three-dimensional communicability space. The positions
of each node correspond to the position vectors $x_{i}$. To enhance clarity, we have associated the first component of the position vectors to the z-axis. Thus, the sign of the $z$-component of each $x_{i}$ informs about
the faction to which each node belongs. That is, the $z=0$ plane separates the nodes into two opposing factions in all cases.
In the unsigned graph (top left), all nodes have a positive $z$ component,
meaning that they all belong to the same faction. Moreover, the value
of the $z$-component could be understood as a measure of the ``polarization''
of the node, or equivalently, its ``commitment'' to the faction. In
the balanced graph (bottom left), the node with two negative edges has
a positive $z$ component, indicating that this node forms a separate
faction. Thus, the faction structure determined by the $z$ component
coincides exactly with the balanced partition for the balanced graph.
This behavior is guaranteed by claim six of Proposition \ref{position_vectors_proposition}. For the unbalanced graphs, the $z$ component also returns a sensible partition. In the top right graph,
the two nodes connected by a negative edge are assigned to different
factions. However, the unbalance in the graph causes the node with a negative z-component
to have a small magnitude of the $z$-component; i.e., a small ``polarization''.
Finally, in the bottom right panel, all nodes are equidistant from each other,
reflecting that they are automorphically equivalent. The
communicability position vectors partition the graph so that one
of the nodes has a positive $z$ component, another has a negative
one, and the last one has a zero component. Thus, the method successfully
detects that the third node cannot belong to either of the two factions.

Because we have previously seen that the magnitude of the communicability position vectors --i.e., the terms $\Gamma_{ii}$ and $\Gamma_{jj}$-- introduces a bias in any communicability distance-based partitioning algorithm into factions,
we will investigate instead the angle spanned by these vectors.

\subsection{Signed communicability angle}

Using the dot product together with Prop. \ref{position_vectors_proposition}, the angle between the position vectors of any two nodes in a signed graph can be expressed in terms of the signed communicability: 
\begin{align}
\cos(\theta_{ij})=\frac{x_{i}\cdot x_{j}}{{||x_{i}}||||x_{j}||}=\frac{\Gamma_{ij}}{\sqrt{\Gamma_{ii}{\Gamma_{jj}}}},
\end{align}
which prompts us to make the following definition.
\begin{definition}
Let $\Sigma$ be a signed graph with $n$ nodes and communicability
matrix $\Gamma$. Let $x_{i}$ be the position vector defined by $x_{i}(k)=\psi_{k}(i)e^{\lambda_{k}/2}$.
Then, the \textit{communicability angle} $\theta_{ij}$ between nodes
$i$ and $j$ is defined as the angle between the position vectors
$x_{i}$ and $x_{j}$. Mathematically, it is given by the following
formula: 
\begin{align} \label{eq:comm_angle}
\theta_{ij}=\cos^{-1}\frac{\Gamma_{ij}}{\sqrt{\Gamma_{ii}\Gamma_{jj}}}.
\end{align}
\end{definition}
We then have the following:
\begin{claim}
The signed communicability angle is bounded as $0^{\circ}\leq\theta_{ij}\leq180^{\circ}$,
where both bounds are attained in balanced complete graphs with a sufficiently
large number of nodes. The lower bound is attained for pairs of nodes in the same
balanced partition and the upper one for pairs in different
balanced partitions.
\end{claim}

\begin{proof}
{
Let $\Sigma$ be a signed complete graph with indicator function $I$. Then, using Corollary \ref{corollary_2} and $\Gamma_{ii}=\dfrac{e^{n-1}+\left(n-1\right)e^{-1}}{n}$, we get

    \begin{equation}
    \theta_{ij}=\cos^{-1}\frac{\Gamma_{ij}}{\sqrt{\Gamma_{ii}\Gamma_{jj}}}= \begin{cases}
    \cos^{-1}\left(\frac{e^{n}-1}{e^{n}+n-1}\right) \quad \textrm{if} \quad I(i)=I(j), \\
    \cos^{-1}\left(\frac{1-e^{n}}{e^{n}+n-1}\right)  \quad \textrm{if} \quad  I(i)\neq I(j),
    \end{cases}
    \end{equation}
which for $n$ sufficiently large tends to either $\theta_{ij}= 0$, when $I(i)=I(j)$, or $\theta_{ij}= 180^{\circ}$, when $I(i)\neq I(j)$.
}
\end{proof}
\begin{remark}
The previous result illustrates the large separability that the communicability
angle achieves between the nodes in a balanced graph: in complete balanced graphs, it places the nodes in the two antipodes of the hypersphere.
\end{remark}

As can be seen in {Fig. \ref{counterexample}, panels (c)-(d)}, the signed communicability
angle produces a correct partition of the nodes in the graph. That is, the communicability
angle between any pair of vertices in the same partition is small, and the
angle between the vertices belonging to different factions is large. Indeed,
in the following result, we prove that this is always the case when we use the signed communicability angle.

\begin{theorem}
\label{angle_clustering_theorem} Let $\Sigma$ be a balanced graph.
Then, the communicability angle between nodes belonging to the same
balanced faction is always smaller than the communicability angle
between nodes belonging to different balanced factions. Mathematically,
if $I(i)=I(j)$ but $I(i)\neq I(k)$, then $\theta_{ij}<\theta_{ik}$,
where $I$ is the indicator function defined in Section \ref{sec:definitions}. 
\end{theorem}

\begin{proof}
By Theorem \ref{comm_balance_theorem}, we know that, since $i$ and
$j$ belong to the same faction, $\Gamma_{ij}>0$. Given that $\Gamma_{ii},\Gamma_{jj}$
are always positive (Proposition \ref{position_vectors_proposition}, claim 1), it follows
that $\cos\theta_{ij}>0$, so $\theta_{ij}<\pi/2$. On the other hand,
since $i$ and $k$ belong to different factions, we have that $\cos(\theta_{ik})<0$
and thus $\theta_{ik}>\pi/2$. Combining both inequalities, we find
that $\theta_{ij}<\theta_{ik}$. 
\end{proof}

{Another important property of the communicability angle is the following:
\begin{proposition} \label{prop:correlation}
    The matrix $\rho$ with entries $\rho_{ij}= \cos(\theta_{ij})$ is a correlation matrix.
\end{proposition}
\begin{proof}
    Recall that a correlation matrix is a symmetric, positive definite matrix with diagonal elements equal to one. Clearly, $\rho_{ij}$ is symmetric. Additionally, $\rho_{ii} = \frac{\Gamma_{ii}}{\sqrt{\Gamma_{ii}\Gamma_{ii}}}=1$. Finally, since $\Gamma$ is positive definite, it follows trivially that $\rho$ must be positive definite too.
\end{proof}
This last proposition implies that all data science algorithms developed for correlation matrices can be safely applied to the cosine of the communicability angle.  }

We now proceed to modify the signed communicability distance $\left|\left|x_{i}-x_{j}\right|\right|^{2}$
to generate a Euclidean distance capable of correctly classifying 
the nodes in a signed graph into its balanced partition. To do so,
we normalize the position vectors $x_{i}$ and $x_{j}$
in the definition of the signed communicability
distance:

\begin{equation} \label{eq:comm_angle_euclidean}
\left|\left|\frac{x_{i}}{||x_{i}||}-\frac{x_{j}}{||x_{j}||}\right|\right|^{2}=\frac{x_{i}\cdot x_{i}}{||x_{i}||^{2}}+\frac{x_{j}\cdot x_{j}}{||x_{j}||^{2}}-2\frac{x_{i}\cdot x_{j}}{||x_{i}||\ ||x_{j}||} :=  d_{\theta}(i,j).
\end{equation}

In this way, we have proved the following.
\begin{proposition}
The quantity $d_{\theta}(i,j)=2-2\cos(\theta_{ij})$ is a squared
Euclidean distance between two nodes $i$ and $j$ of $\Sigma$. 
\end{proposition}

In the following part of
this work, we proceed to use these metrics to solve some common data science problems, like the partitioning of unbalanced networks into factions.

\section{Applications}

In previous sections, we have seen that the communicability provides
several well-defined measures for signed networks: among others, the communicability
distance and the communicability angle. These measures can be interpreted
as the separation of nodes in an ``alliance space''. This alliance
space tends to group pairs of vertices forming alliances close together, separated from those with whom they maintain antagonistic relations.
Moreover, simple transformations, like taking the cosine of the communicability
angle, can transform these measures into similarity metrics. In this
section, we will explore how these measures can be used in connection
with simple statistical and machine-learning techniques to gain novel
insights into the structure of empirical signed networks.

{

The framework we will use to analyze empirical signed networks is inspired by \cite{pereda2019visualization} and based on the communicability angle.
Specifically, we leverage Theorem \ref{angle_clustering_theorem}, which showed that the communicability angle correctly identifies 
the balanced factions of a structurally balanced network. Because of this, $\theta_{ij}$ is ideally suited to analyze structures such as factions present in real-world data. Moreover, the communicability angle 
is a valid dissimilarity metric. This means that we can use it as input for 
unsupervised algorithms like agglomerative clustering or multidimensional scaling.
Building on this, we propose the following framework for analyzing signed networks,
consisting of three main steps:

\begin{enumerate}
\item Compute the communicability angle matrix from the signed adjacency matrix $A$.
\item Use a dimensionality reduction algorithm to embed the network into a low-dimensional space, using the communicability angle as the input dissimilarity matrix\footnote{{Some dimensionality reduction algorithms, like PCA, use a correlation matrix as input. In those cases, one should instead use the cosine of the communicability angle, which is a well-defined correlation matrix (Prop. \ref{prop:correlation}).}}.
\item Perform a clustering algorithm in the low-dimensional embedding to identify network factions.
\end{enumerate}

This framework mirrors the procedure followed in spectral clustering, where the data is embedded 
into a low-dimensional space using Laplacian eigenvectors, and then clustered using K-means or any other
clustering algorithm
 \cite{vonluxburgTutorialSpectralClustering2007}. The primary difference in our approach lies in using the communicability angle matrix instead of a conventional affinity matrix. This enables the use of embedding methods that do not necessarily rely on the Laplacian eigenvectors, overcoming the limitations of spectral clustering. Moreover, the communicability angle aggregates local and non-local information
about the network structure, which helps avoid local minima in the optimization steps within the embedding and clustering algorithms. 

A key advantage of the proposed algorithm is its flexibility: each of the three steps is compatible with different metrics and algorithms, thereby allowing to exploit the strengths of each one depending on the context. For instance, to perform dimensionality reduction, one can use PCA, t-SNE, multidimensional scaling, etc {(see for instance \cite{izenman2008modern} for an overview of dimensionality reduction algorithms)}. Similarly, for the clustering step,
one can choose between K-means, DBSCAN, hierarchical clustering, spectral clustering, etc. Finally, it is not even required to use the communicability angle as input (although it is recommended). One could use the communicability distance, the resistance distance, or any other well-defined dissimilarity metric.
}

{In the remainder of this work},
we will focus on applications {of this framework} in (i) network visualization, (ii) detection
of factions of allies, (iii) determining a hierarchy of alliances,
(iv) dimensionality reduction, and (v) quantification of node polarization.
We will exemplify these techniques by using four distinct datasets:
firstly, a network of alliances and enmities between tribes in New
Guinea (Gahuku-Gama dataset); secondly, a temporal network of international relations throughout
the XX century (IR dataset); thirdly, a temporal network of voting
correlations in the European Parliament (European Parliament dataset); {and finally, the gene regulatory network of \textit{E. coli}}. {For the sake of simplicity, we will always choose 
the communicability angle matrix as the dissimilarity metric, metric multidimensional 
scaling (MDS) as the dimensionality reduction algorithm, and K-means or 
hierarchical clustering as the clustering algorithms. 
The embedding dimension in
the MDS algorithm
is always equal to two, so we can visualize the results in a 2D plot.
To find the optimal
number of clusters, we will use silhouette index maximization; this is, we compute the average silhouette index $\bar s$ for a range of cluster numbers and select the number that maximizes $\bar s$.
When using hierarchical clustering, we choose an agglomerative procedure with Ward or average linkage \cite{scikit-learn}.
}

\subsection{Gahuku-Gama tribes}

\begin{figure}[h]
\centering 
\includegraphics[width=0.7\textwidth]{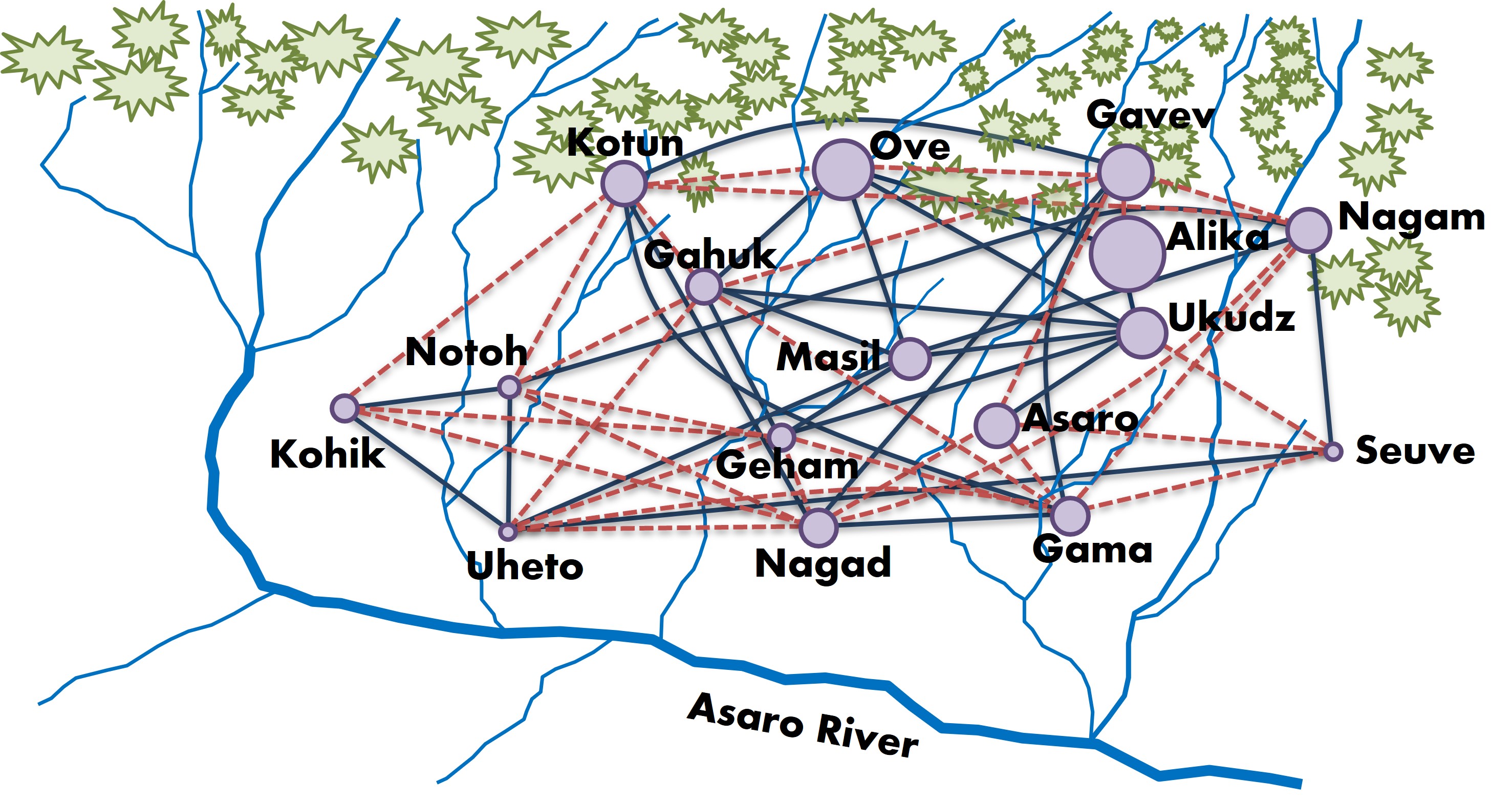} \caption{{Sketch of the Gahuku-Gama network. The nodes are located in the geographical
region that the corresponding tribe inhabits.}}
\label{fig:Gahuku-Gama}
\end{figure}

To warm up, we start by considering the interactions between 15 indigenous
tribes located in the Eastern Highlands of the island of Papua New
Guinea \cite{read1954cultures}. The analysis of the ``Gahuku-Gama'' tribal system has
previously been studied by using traditional anthropological techniques
\cite{read1954cultures}, and the current analysis exemplifies
the use of data analysis based on signed networks in anthropology.
One notable aspect of the Gahuku-Gama society is the continuous state
of conflict between the different tribes. Violent conflicts arise
every dry season as a consequence of land disputes and accusations
of sorcery. Moreover, the tribes often create alliances to gain an
advantage against their enemy tribes, forming a complex web of alliances
and hostilities. Consequently, the Gahuku-Gama network offers an invaluable
testing ground for mathematical theories regarding signed interactions. The
geographical distribution of the Gahuku-Gama tribes is shown in Fig. {\ref{fig:Gahuku-Gama}}. Notice that, based only on geographical location,
it is difficult to infer which tribes form a cohesive faction.

We begin by showing how the communicability angle matrix can provide
a meaningful visualization of the Gahuku-Gama network and how this
can be leveraged to infer antagonistic factions in the network. To
do so, we will employ the multidimensional scaling algorithm. 
The resulting embedding
is shown in Fig. \ref{fig:Gahuku_Gama_results}. Contrary to the
geographical visualization of the network, our communicability-based
embedding clearly depicts the faction structure of the network. Firstly,
we can see that there is a faction composed of four tribes: GAVEV,
GAMA, KOTUN, and NAGAD, whose defining characteristic is that they
have a common enemy: the NAGAM tribe. The rest of the tribes are split
into two additional factions, with the NAGAM tribe acting as a bridge
between them. To confirm these qualitative results, we have run a
$K$-means algorithm on the two-dimensional embedding coordinates
to detect the cluster structure of the network. The algorithm confirmed
that the three-faction structure is indeed the optimal one according
to the minimization of the Silhouette index \cite{rousseeuw1987silhouettes}.
The clusters returned by the $K$-means algorithm are indicated through
the node colors of Fig. \ref{fig:Gahuku_Gama_results}(a).

Additionally, the communicability angle also provides a simple method
for establishing a hierarchy of alliances between the different tribes.
By this, we mean a hierarchical structure that provides a ranking
of alliances between the pairs of tribes, from the strongest alliance
to the weakest one. To achieve this goal, we employ
a hierarchical clustering algorithm using the Euclidean distance on the 2D embedding coordinates.
This procedure results in the dendrogram depicted in Fig. \ref{fig:Gahuku_Gama_results}(b), which confirms the three-faction structure that was found
both visually and through a clustering algorithm. Moreover, it provides
a more nuanced analysis of the emerging alliances between the system
of Gahuku-Gama tribes. For example, it ranks the ASARO-OVE-MASIL and
NOTOH-UHETO as the two strongest groups of allies. Surprisingly, OVE
and ASARO are not connected through an alliance. Therefore, its remarkably
strong alliance is an emerging feature that results from their common
allies and enemies. This example clearly shows how effective alliances
established through common allies and enmities can have a stronger
effect than pairwise alliances. Consequently, local measures such
as the classical centrality measures can fail when assessing the hierarchy
of alliances, whereas measures that account for the global structure
of the signed network, such as the communicability, will provide more
accurate results.

\begin{figure}
\subfloat[]{
    \fcolorbox{black}{white}{
    \includegraphics[width=0.46\textwidth]{figures/Gahuku_Gama_MDS} 
    }}
\subfloat[]{
    \fcolorbox{black}{white}{
    \includegraphics[width=0.46\textwidth]{figures/Gahuku_Gama_dendrogram} 
    }}
\par
\centering{}
\caption{Representation of the partitions of the Gahuku-Gama network, using multidimensional scaling into a two-dimensional space (panel (a)) followed by a K-means clustering algorithm (panel (a)) and hierarchical clustering (panel (b)).}
\label{fig:Gahuku_Gama_results} 
\end{figure}

\subsection{International relations in wartime}

\begin{figure}[h]
\centering{}\includegraphics[width=0.8\textwidth]{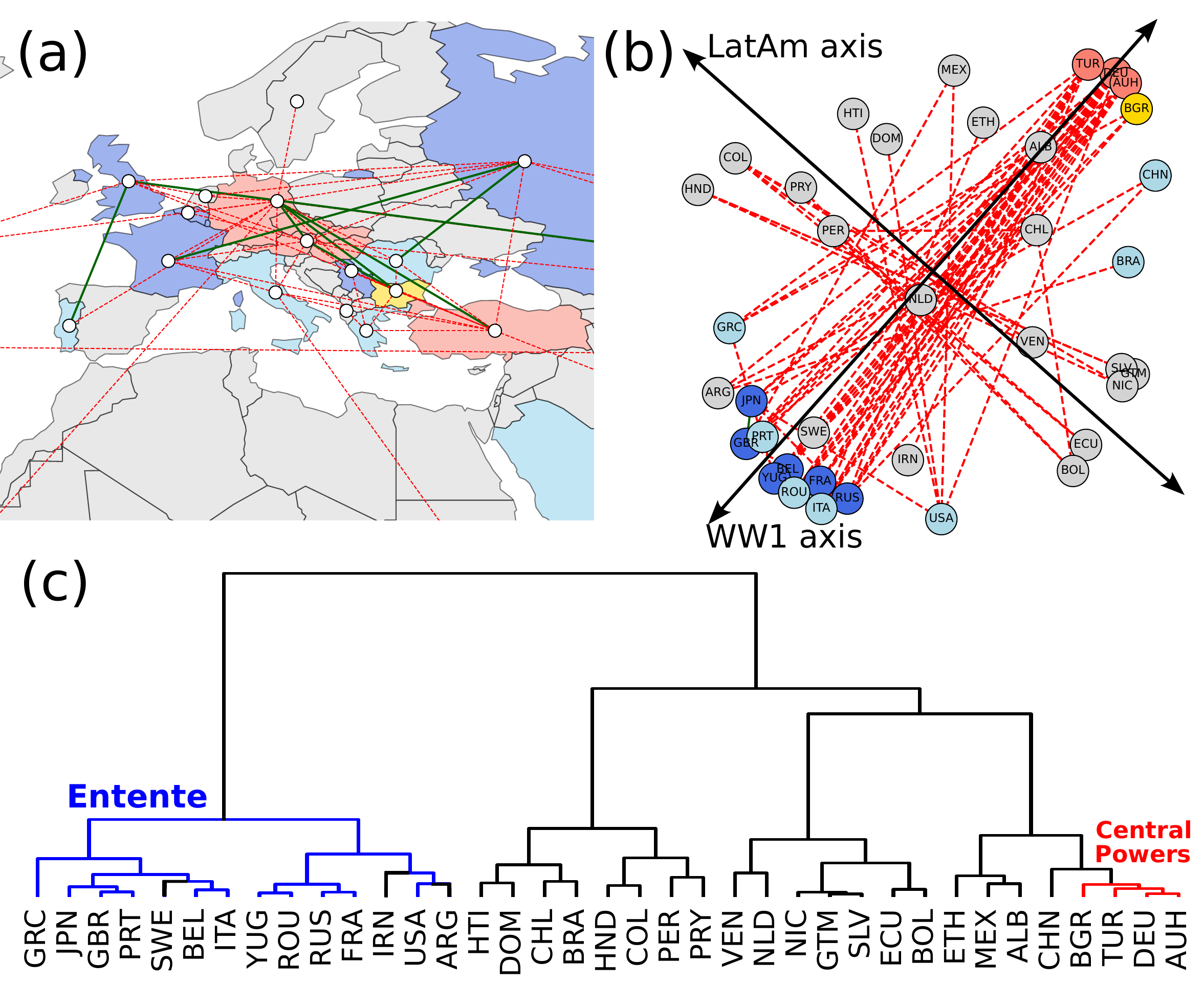} \caption{International relations in 1914 represented as a signed graph (a), 2D embedding in the communicability angle space (b), and corresponding dendrogram using ward linkage (c).}
\label{fig:1914}
\end{figure}

The second analysis focuses on the use of signed networks for data
analysis of international relations (IR). {We use the networks of alliances \cite{gibler2008international}, military interstate disputes \cite{palmer2022mid5}, and strategic rivalries \cite{colaresi2008strategic}, obtained from the Correlates of War database (see \cite{diaz-diaz2023} for details regarding the construction of the network)}. Let us focus on the First and Second World Wars, where the web of IR defines a clear faction structure, which far from being static, was changing over time due
to the incorporation of countries into the war effort, cease of hostilities,
or shifting alliances. We start by analyzing international relations
during 1914, corresponding to the beginning of World War I (WWI).
The geographical and communicability embeddings of the network are
shown in Figs. \ref{fig:1914}(a) and (b). Dark blue and red colors
denote the original members of the Entente and Central Powers, respectively.
Again, the geographical embedding does not give a clear picture of
the faction structure emerging in 1914, whereas the communicability
embedding clearly clusters countries according to the faction they
joined during WWI. In particular, all the original members of the
Entente in 1914 (Serbia, Russia, France, UK, Belgium, and Japan) are
located in the bottom left corner, while the original members of the
Central Powers (Germany, Austria-Hungary, and the Ottoman Empire) are
located in the upper right corner. More surprisingly, each cluster
of the embedded network contains most of the countries that would
join the corresponding faction after 1914 (highlighted with light
blue and gold colors). In other words, despite the formal neutrality
that these countries declared in 1914, the communicability classifies
them as aligned with a given faction, which turns out to be the faction
they eventually join in future years. The associated hierarchical
dendrogram (Fig. \ref{fig:1914}(c)) confirms these results by grouping
most members of the Entente in the left branch and all Central Powers
in the right branch.

The two-dimensional embedding discerns a second component present
during this period, which is orthogonal to the one formed by the powers
at war during WWI. This second axis is composed of Latin American
countries, which are well-known to remain neutral during this conflict.
Historians have recognized that only about twenty countries remained
neutral during WWI, mainly Latin American and Northern
European ones \cite{rinke2017latin}. The 2D embedding generated with the communicability thus reveals the existence of two
``principal axes'' of conflict, in the same way that other statistical
techniques, such as the principal component analysis (PCA), discern
the principal components of a multivariate dataset.

\begin{figure}[h]
\centering{}\includegraphics[width=1\linewidth]{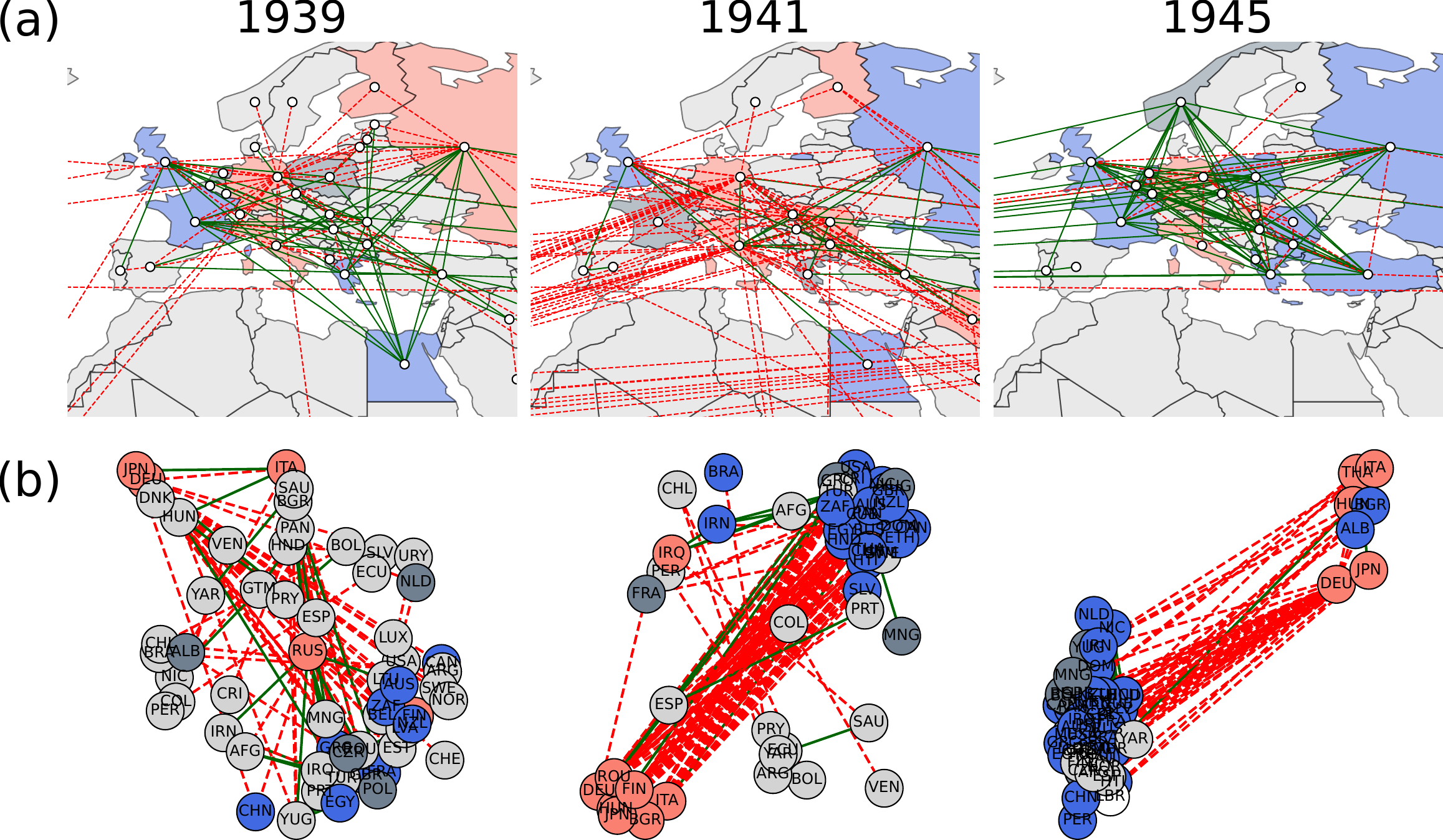} \caption{(a) Network of international relations in Europe during World War
II, in 1939 (left column), 1941 (middle column) and 1945 (right column).
Blue color indicates members of the Allies, red color indicates members
of the Axis, and dark gray colors indicate occupied countries. (b)
Communicability embedding corresponding to each of the networks.
Same color code as (a).}
\label{fig:Embedding_WW2}
\end{figure}

We now move to the analysis of World War II (WWII), where a more complicated
pattern emerges. Instead of a clear division of factions since the
beginning of the conflict, one can observe an initially unpolarized
network with unclear factions, which becomes more polarized as the
war progresses (see Fig. \ref{fig:Embedding_WW2}). One key player
in this dynamic is the Union of Soviet Socialist Republics (USSR, iso code RUS). Indeed, in 1939, it occupied
the most central position on the embedding\footnote{{A comprehensive introduction to the complex web of relations during World War II can be found in \cite{beevor2012second}.}}. This central position
reflects the unorthodox diplomatic choices of the USSR during this
year: despite the markedly anti-communist ideology of Nazism, Germany
and the USSR signed the Ribbentrop-Molotov non-aggression
pact in 1939, which secretly agreed to partition Central and Eastern
Europe between them. Thus, the node corresponding to the USSR tends
to mark its distance from both the fascist states and Poland's allies,
resulting in the observed central position. In 1941, Nazi Germany
turned against the USSR and began Operation Barbarossa, resolving
the ``tension'' in the network structure and eventually causing
the emergence of a clear faction structure. The factions emerging
in 1941 correspond to the Axis Powers and the Allies. Note the unusual
positioning of France in the central region, which reflects the Nazi
occupation of Northern France and the creation of the collaborationist
regime of Vichy France. The faction structure of the network increases
further until 1945. In this final year of war, virtually every country
had joined one of the two factions, having a clear two-faction structure.
The unusual positioning of Albania and Bulgaria reflects that, even
though they were nominally Allies in 1945, they had been part of the
Axis power until the previous year: Bulgaria was a German ally until
a communist coup d'etat in 1944 overthrew its Government, while Albania
was liberated from German occupation in November 1944. This example
illustrates the potentialities of the use of communicability functions
of signed networks in the data analysis of historic events, which
may be a point in the direction of a more quantitative analysis of
history.

\subsection{Voting system at the European Parliament}
\begin{figure}
\centering{}\includegraphics[width=1\textwidth]{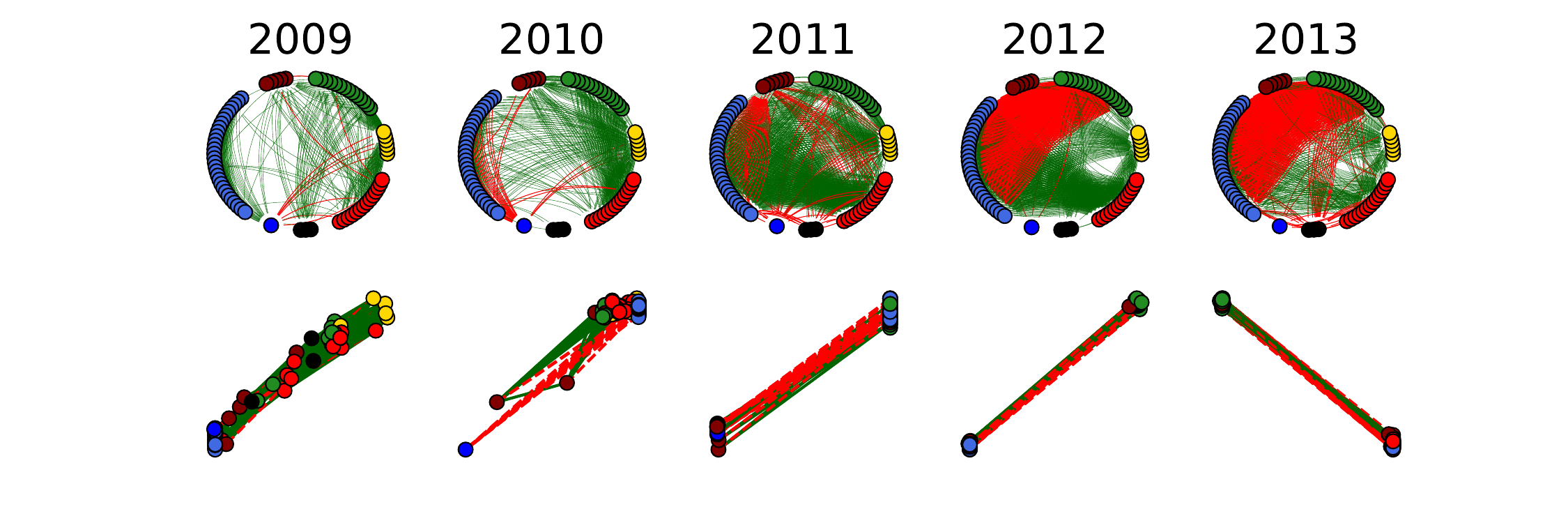}
\caption{Illustration of the signed networks of the voting system of French
members of the European Parliament from 2009 to 2013 in topics related to agriculture.
In the top line, we show the networks drawn in a circle, coloring
the representatives of different groups with different colors. In the
bottom line, we illustrate the 2-dimensional projection of the same
networks based on the communicability geometry.}
\label{fig:EU_Parliament}
\end{figure}

\begin{figure}
\centering{}\includegraphics[width=1\textwidth]{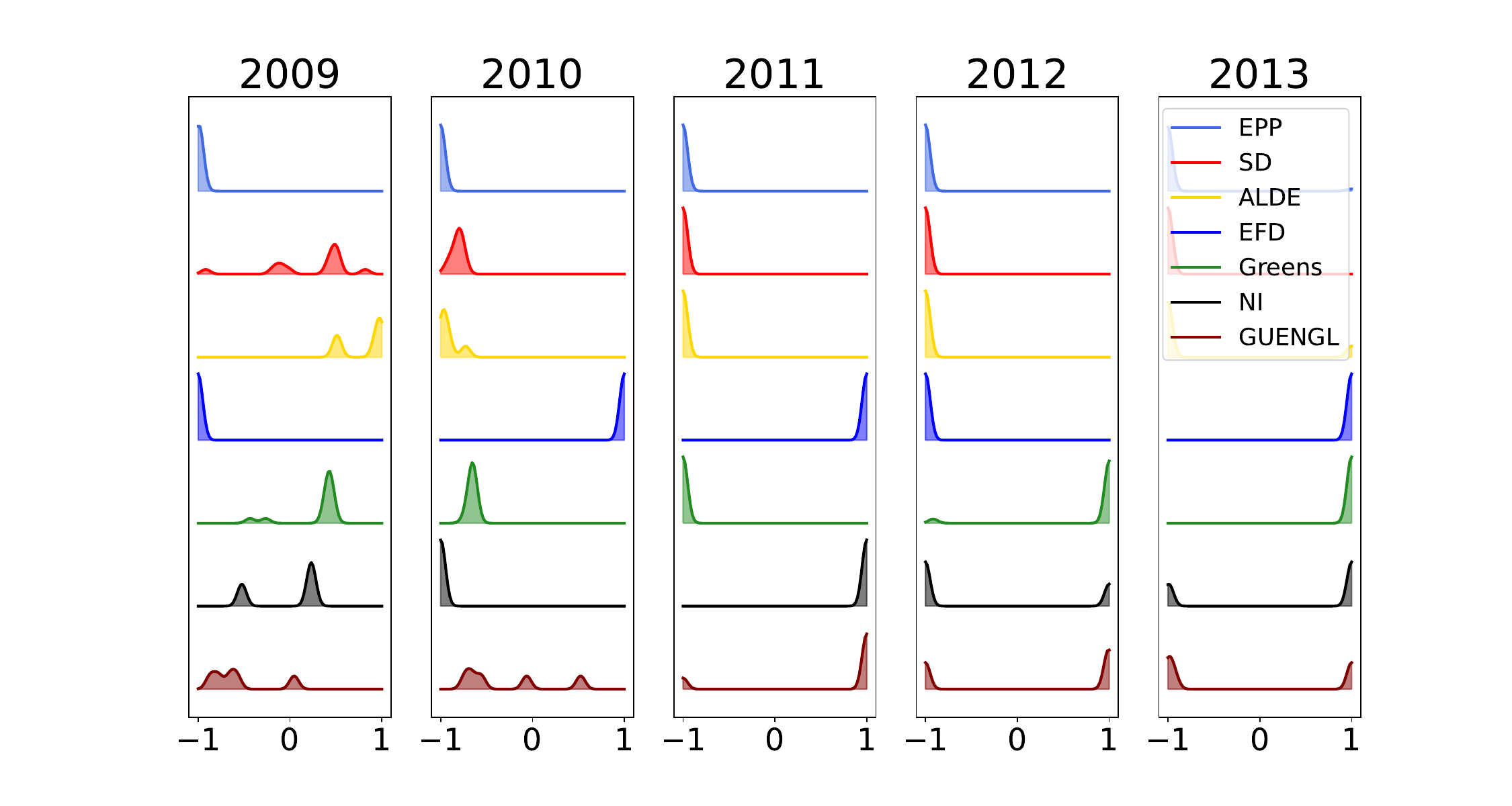}
\caption{{Latent polarization distribution of each political group for each
network analyzed in the time range 2009-2013. The latent polarization
is quantified as the first component of a PCA of the embedding shown
in Fig. \ref{fig:EU_Parliament}. The PDFs were smoothed using kernel
density estimation with Gaussian kernels and bandwidth $h=0.05$.}}
\label{fig:EU_Parliament_hist}
\end{figure}

The third case study on the applications of communicability-based
data analysis of signed networks corresponds to its application to
political sciences. One way of representing political systems using signed networks is  by modeling political actors as nodes and their affinities or hostilities (recorded by their frequency of common/opposed votes) as signed links. Such ``correlation'' in their voting pattern
can be quantified by means of the Pearson correlation coefficient
(other conventions to define the edges can be chosen, see \cite{masuda2023correlation}).
A frequently appearing feature of these systems is the emergence of
\textit{polarization}, where two political factions with completely
anticorrelated voting patterns emerge. 

Here, we analyze the temporally evolving network constructed from the voting similarities among Members of the European Parliament (MEPs). {We employ the data and network construction protocol from} \cite{arinik2017signed, arinik2020multiple}. This protocol entails three steps: (1) Collection
of voting records for each MEP across roll-call votes during the 7th
term of the European Parliament (2009-2014), sourced from the ``It's
Your Parliament'' website\footnote{http://www.itsyourparliament.eu/}. Our analysis focuses on policy domains
about agriculture, and specifically considers MEPs representing
France. This selection is motivated by France's significant stake
in the Common Agricultural Policy (CAP) during this term, given the
importance of agriculture to its economy. (2) Calculation of the similarity
between all pairs of MEPs within a given political year. This similarity
measure is taken to be proportional to the number of times two MEPs
voted equally (for example, both voting against a given law) minus
the number of times they did not vote equally. To ensure that the
similarity lies between -1 and +1, it is divided by the total number
of roll-calls of the political year. (3) Application of a thresholding
procedure to eliminate weak similarities, thereby converting the dense
similarity matrix into a sparse signed network. Further details on
this thresholding method can be found in \cite{arinik2017signed}. The nodes
of the network belong to one of seven political parties: EPP (center-right,
blue color), SD (center-left, red color), ALDE (center, gold color),
EFD (right, dark blue color) Greens (environmentalist-left, green
color), GUE-NGL (far-left, dark red color), and NI (independent, often
far-right, black color). We begin the analysis by embedding the network
in a plane, using {our proposed framework}. In Fig. \ref{fig:EU_Parliament}, we compare the embedding
with a circular layout with the nodes clustered based on their political
group. Surprisingly, the political group does not reveal a clear partition
of the network, due to the existence of several negative correlations
within each Parliamentary group and positive correlations between
different ones. In contrast, the communicability embedding unravels
several hidden features of this political system. On the other hand, for the communicability embedding,
it is remarkable that although the projection of the data is on a
two-dimensional space, the vertices tend to align along a straight
line, particularly as time evolves towards 2013. This result contrasts
with the latent ideological spaces of other political entities, which
seem to have larger dimensions (see e.g. \cite{martin-gutierrezMultipolarSocialSystems2023}). Furthermore, the distribution of the points along
this straight line changes notably from 2009 to 2013. In the first
period, we find an almost homogeneous distribution of the points along
the line, while in 2013 we find that every point is at each of the
two extremes. This finding clearly points out an increasing
tendency in the polarization of French representatives in the European Parliament. 

Because the optimal embedding aligns the nodes
along a one-dimensional axis, we can utilize the nodes' positions
on this axis as a measure of its latent polarization. Formally, this
corresponds to performing a PCA of the embedding coordinates and taking
the first component of every node as its polarization. Moreover, we
normalize the components so that their values range between
-1 and +1, so that we can compare node's polarization for different
years. The probability distributions of these latent polarizations,
grouped by political affiliation, are illustrated in {Fig. \ref{fig:EU_Parliament_hist}}.
The probability density functions once again show a growing polarization over time. Interestingly, the representatives of the three largest
parties (EPP-SD-ALDE) are mostly grouped within the same factions
while the second faction has a volatile composition, with representatives
from EFD, the Greens, NI, and GUE-NGL coming in and out. Let us
focus on the Greens, which evolve from a weakly polarized position
in 2009 to a large polarization in 2013, where they aligned with EFD,
which is traditionally in the opposite faction. In 2013, the European Parliament had
to vote for a new Common Agricultural Policy (CAP) {\cite{ReformCommonAgricultural}}. This is an extremely
important topic in France in which agriculture plays a fundamental
role in the economy. In this case, the Greens qualified the new proposal
of CAP as unjust from an environmental point of view and voted against
it, thus aligning with the conservatives of the EFD, who were also opposing the
law. These results show that the European Parliament dynamics are
highly complex and strongly defy the right-left traditional narrative.
More importantly, it serves as a clear illustration of the potentialities
of the communicability-based analysis of signed networks in political
data analysis.

\subsection{Gene regulatory network of \textit{E. coli}}

{
As a final case study, we will consider the gene regulatory network of \textit{Escherichia coli} (\textit{E. coli}). This network can be naturally modeled as a signed graph, where nodes represent genes or transcription factors, and edges capture regulatory interactions. The sign of an edge encodes the nature of regulation: positive edges correspond to activation, while negative edges represent inhibition. Unlike social networks, where negative links typically indicate antagonism and are perceived as undesirable, inhibitory interactions in biological systems play a key role in maintaining homeostasis. Despite this difference in interpretation, concepts from structural balance theory remain relevant. Specifically, balanced motifs promote stability and irreversibility, whereas unbalanced motifs may lead to oscillatory or switch-like behaviors \cite{alon2019introduction}. Analyzing the structure of the signed interactions is thus essential for understanding the functional organization of transcription networks.

For this analysis, we use the dataset from \cite{aref2019balance}, based on the RegulonDB database \cite{salgado2006regulondb}. The analysis is restricted to the main connected component, which comprises 1,376 nodes and 3,150 edges, of which 1,848 (59\%) are positive and 1,302 (41\%) are negative.  The network has been symmetrized, with edges displaying ambiguous sign patterns removed. To uncover the latent structure of this network, we compute the communicability angle matrix. Remarkably, the matrix displays a clear block structure, characterized by positive diagonal blocks and negative off-diagonal blocks (Fig. \ref{fig:ecoli}, left panel). This block structure is further validated by the dendrogram shown in Fig. \ref{fig:ecoli} (right panel), generated using an agglomerative clustering algorithm with average linkage}\footnote{{Since the block structure is already evident in the high-dimensional space, we chose to omit the framework's dimensionality reduction step.}}.
{The dendrogram reveals two major blocks, with the lower block subdivided into two smaller modules. This three-block structure suggests the presence of distinct regulatory modules or factions within the network, each characterized by internally coherent patterns of activation and inhibition. Furthermore, the structure appears hierarchical, with larger blocks containing smaller nested modules. Biologically, this partitioning may reflect distinct regulatory circuits, such as modules governing metabolism, stress responses, or cell cycle processes. In summary, our communicability-based analysis also provides a powerful approach in biological contexts, where it could identify and characterize the regulatory modules that govern functions like development, differentiation, and response to the environment.

\begin{figure}
\centering{}\includegraphics[width=1\textwidth]{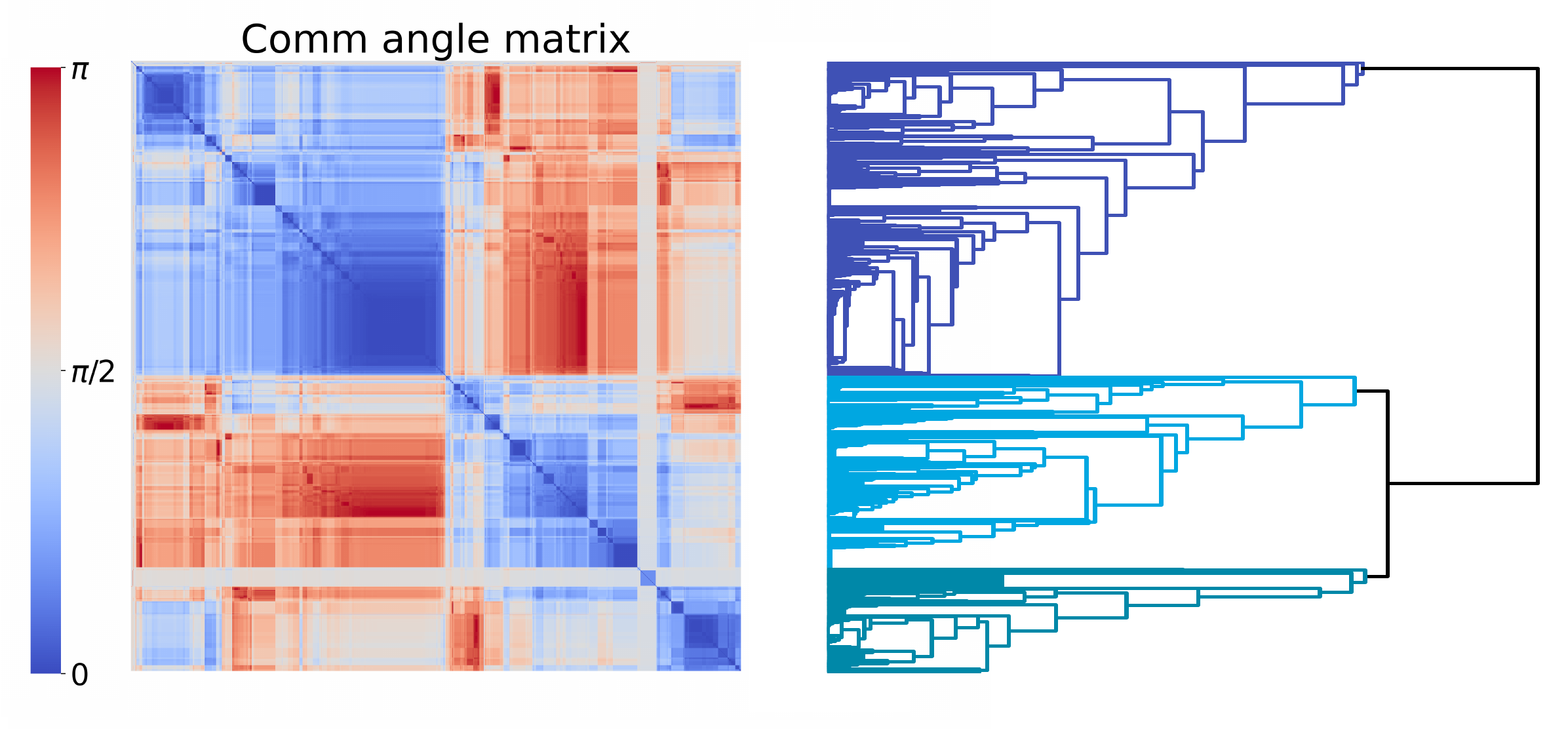}
\caption{{Left panel: communicability angle matrix of the \textit{E. coli} transcriptional network. Blue shades represent communicability angles smaller than $\pi/2$ (effective allies) while red colors correspond to angles greater than $\pi/2$ (effective enemies). Right panel: dendrogram generated using agglomerative clustering with average linkage, based on the communicability angle matrix. Clusters are highlighted in distinct colors for clarity. }}
\label{fig:ecoli}
\end{figure}

}

\section{Conclusions and future work}

{

This study introduces a novel framework for analyzing signed networks through the concept of communicability geometry, providing robust mathematical tools to address key challenges in the data analysis of signed graphs.  Our main contribution is the introduction of the signed communicability distance and angle. These metrics, grounded in the exponential of the signed adjacency matrix, satisfy the axioms of a distance while incorporating information from both positive and negative edges. Second, we showed that the communicability induces a hyperspherical embedding for signed networks, offering a geometric representation of the relationships between nodes in a signed graph. Importantly, this geometric representation opens the door to the use of common data analysis tools, like dimensionality reduction. This leads us to the third contribution of this work, which is the proposal of a flexible framework to analyze signed networks through the lens of communicability. Using the proposed metrics, we address diverse problems in data analysis of signed networks, including network partitioning, dimensionality reduction, hierarchy detection, and polarization quantification, in a unified and mathematically coherent way.

The key strengths of the proposed framework are its mathematical rigor, flexibility, and applicability across different domains. Unlike previous approaches, ours respects key mathematical properties, like the distance axioms. Moreover, it guarantees the identification of the balanced factions when the graph is balanced (Theorem \ref{angle_clustering_theorem}). Another notable advantage of the communicability framework is that it provides not only well-defined distance matrices but also correlation matrices (Proposition \ref{prop:correlation}), enabling the application of techniques designed for correlation-based analysis. Furthermore, the framework integrates with many machine learning and data analysis algorithms, such as multidimensional scaling, hierarchical clustering, PCA, DBSCAN, Planar Maximally Filtered Graphs, UMAP, and many others. This flexibility extends its applicability to a wide range of real-world scenarios, from gene regulatory networks to international relations, to name a few. 

This work is not without limitations. First, like any model, the representation of a complex system through a signed network abstracts away certain nuanced aspects, such as node attributes, neutral edges, and multidimensional relationships. Additionally, our analysis is restricted to simple signed networks, leaving out more complex network types, such as directed, weighted, or multilayer networks, which are often encountered in real-world data. Another limitation lies in scalability; the direct computation of the matrix exponential could be computationally expensive for large datasets. Fortunately, a variety of methods exist to address this issue \cite{moler2003nineteen}. For example, Krylov subspace methods, such as the Lanczos or Arnoldi processes, allow for efficient approximation of the matrix exponential by focusing on its action on a vector rather than explicitly forming the entire matrix. Similarly, algorithms based on truncated series expansions, such as Taylor or Chebyshev expansions, can provide computationally efficient approximations while preserving sparsity.

This work primarily establishes a solid mathematical foundation for analyzing signed networks, laying the groundwork for diverse data analysis procedures. Future efforts will focus on practical aspects of implementing and optimizing the framework. Specifically, a detailed evaluation of the optimal combination of distance metrics, dimensionality reduction techniques, and clustering algorithms is necessary to improve the framework's accuracy and robustness in real-world applications. Another important direction is optimizing the framework's implementation to ensure scalability for large datasets. Additionally, future work will benchmark the framework against existing faction detection techniques to clarify its advantages and limitations. Finally, we intend to extend the framework to classification and prediction tasks. These efforts will ensure the practical applicability of the framework while maintaining its strong theoretical foundation.
}

\section*{Acknowledgments}
EE thanks Zeev Maoz for sharing data.
FDD and EE acknowledge funding from Ministerio de Ciencia e Innovacion, Agencia Estatal de Investigacion Program for Units of Excellences Maria de Maeztu (CEX2021$-$001164-M$/$10.13039$/$501100011033).  FDD thanks financial support  from MDM$-2017-0711-20-2$ funded by MCIN$/$AEI$/$10.13039$/$501100011033 and by FSE invierte en tu futuro, as well as project APASOS (No. PID2021 $-$122256NB$-$C22). EE also acknowledges funding from project OLGRA (PID2019-107603GB-I00) funded by Spanish Ministry of Science and Innovation.

\bibliography{references}
\bibliographystyle{ieeetr}

\end{document}